\documentclass[a4paper, 11pt, twoside]{article}
\usepackage[nogytt]{mfpaperstuff2}
\usepackage{tikz-3dplot}
\usetikzlibrary{shapes}

\newgeometry{margin=0.421in}
\begin{document}
\newcommand\three[4]{\draw(#1,0)++(120:#2)++(240:#3)node[label={above,fill=white,inner sep=1pt:{$#4$}},fill,circle,inner sep=2pt]{}node[label={below,fill=white,inner sep=1pt:{$#1,#2,#3$}}]{};}
\newcommand\four[5]{\draw(#1,#2)node[label={above,fill=white,inner sep=1pt:{#5}},fill,circle,inner sep=2pt]{}node[label={below,fill=white,inner sep=1pt:{$#1,#2,#3,#4$}}]{};}
\newcommand\five[6]{\draw(#1,#2)node[label={above,fill=white,inner sep=1pt:{#6}},fill,circle,inner sep=2pt]{}node[label={below,fill=white,inner sep=1pt:{$#1,#2,#3,#4,#5$}}]{};}
\newcommand\ssc[2]{\operatorname{SC}_{#1}(#2)}
\newcommand\stc[2]{\operatorname{SD}_{#1}(#2)}
\newcommand\sz[1]{\widebar{#1}}
\newcommand\act[1]{\mathring#1}
\newcommand\zzsz{\{0,\dots,s-1\}}
\newcommand\zztz{\{0,\dots,t-1\}}
\newcommand\spr{$s$-periodic\xspace}
\newcommand\kere[1]{K^{(#1)}}
\newcommand\kerd[1]{L^{(#1)}}
\newenvironment{thmcx}[1]{\begin{thmciting}{\textup{\textbf{\cite{#1}.}}\ }}{\end{thmciting}}
\medmuskip=3mu plus 1mu minus 1mu
\newcommand\idx[2]{\card{#1\,{:}\,#2}}
\newcommand\ci[2]{#1{\circ}#2}
\newcommand\stab[3]{\operatorname{Stab}_{#1,#2}(#3)}
\newcommand\rstb[2]{\operatorname{Stab}_{s,t}^{(#1)}(#2)}
\newcommand\rstbsc[2]{\operatorname{StabSC}_{s,t}^{(#1)}(#2)}
\newcommand\stg[1]{\operatorname{Stab}_G{#1}}
\newcommand\stb[3]{\card{\stab{#1}{#2}{#3}}}
\newcommand\stabsc[3]{\operatorname{StabSC}_{#1,#2}(#3)}
\newcommand\stbsc[3]{\card{\stabsc{#1}{#2}{#3}}}
\newcommand\ltup[2]{\left(\left.#1\ \right|\ \smash{#2}\right)}
\newcommand\rtup[2]{\left(\smash{#1}\ \left|\ #2\right.\right)}
\newcommand\lmst[2]{\left[\left.#1\ \right|\ \smash{#2}\right]}
\newcommand\rmst[2]{\left[\smash{#1}\ \left|\ #2\right.\right]}
\newcommand\plu[2]{#1\boxplus#2}
\newcommand\weyl{\tilde{\mathfrak S}_}
\newcommand\lyew{W_}
\newcommand\wsn{\weyl s^{(n)}}
\newcommand\hyp{\tilde{\mathfrak H}_}
\newcommand\lra\longrightarrow
\newcommand\cornoarg[1]{\mathrm{cor}_{#1}}
\newcommand\weinoarg[1]{\mathrm{wt}_{#1}}
\newcommand\cor[2]{\cornoarg{#1}#2}
\newcommand\wei[2]{\weinoarg{#1}#2}
\newcommand\quo[2]{#2^{(#1)}}
\newcommand\qou[1]{\mathrm{quo}_{#1}}
\newcommand\be[1]{\mathcal{B}^{#1}}
\newcommand\ber[2]{\calb_{#1}^{#2}}
\newcommand\gc[2]{$[#1{:}#2]$-core}
\newcommand\gcc[3]{[$#1{:}#2{:}#3$]-core}
\newcommand\cg[2]{\mathcal{C}_{#1:#2}}
\newcommand\cores[1]{\calc_{#1}}
\newcommand\sccores[1]{\cald_{#1}}
\newcommand\rcores[2]{\cores{#1}^{(#2)}}
\newcommand\rsccores[2]{\sccores{#1}^{(#2)}}
\newcommand\rcn{\rcores s{nt}}
\newcommand\rcnn{\rcores sN}
\newcommand\rdnn{\rsccores sN}
\newcommand\rdn{\rsccores s{nt}}
\newcommand\cnd{(\textasteriskcentered)}
\newcommand\cne{(\textdagger)}
\newcommand\zsz{\bbz/s\bbz}
\newcommand\ztz{\bbz/t\bbz}
\newcommand\es[2]{\Gamma_{#1}#2}
\newcommand\twq[1]{\mathrm Q_s^t#1}
\newcommand\sset[1]{\cals(#1)}
\newcommand\simcores[2]{\cores{#1}\cap\cores{#2}}
\newcommand\simsccores[2]{\sccores{#1}\cap\sccores{#2}}
\newcommand\eqrn[1]{\equiv_{#1}}
\newcommand\neqrn[1]{\nequiv_{#1}}
\newcommand\eqr{\eqrn n}
\newcommand\neqr{\neqrn n}
\newcommand\eqn{\eqrn N}
\newcommand\neqn{\neqrn N}
\newcommand\eqnt{\eqrn{nt}}
\newcommand\neqnt{\neqrn{nt}}
\newcommand\lb[1]{node[midway,fill=white,inner sep=1pt]{\footnotesize$\sz#1$}}
\newcommand\lc[1]{node[pos=.25,fill=white,inner sep=1pt]{\footnotesize$\sz#1$}}
\newcommand\ld[2]{node[pos=#2,fill=white,inner sep=1pt]{\footnotesize$\sz#1$}}
\newcommand{\mathsc}[1]{{\normalfont\textsc{#1}}}
\newcommand\sci{\mathsc i}
\newcommand\scj{\mathsc j}

\pdfpagewidth=12.5cm
\pdfpageheight=8.839cm
\setcounter{page}0
{\footnotesize This is the author's version of a work that was accepted for publication in the Electronic Journal of Combinatorics. Changes resulting from the publishing process, such as peer review, editing, corrections, structural formatting, and other quality control mechanisms may not be reflected in this document. Changes may have been made to this work since it was submitted for publication. A definitive version was subsequently published in\\
\textit{Electronic J. Combin.} \textup{\textbf{23} (2016) \#P4.32.\\
{\scriptsize http://www.combinatorics.org/ojs/index.php/eljc/article/view/v23i4p32/pdf}}}
\newgeometry{margin=1in,includehead,includefoot}

\title{$(s,t)$-cores: a weighted version of Armstrong's conjecture}
\toptitle

\begin{abstract}
The study of core partitions has been very active in recent years, with the study of \emph{$(s,t)$-cores} -- partitions which are both $s$- and $t$-cores -- playing a prominent role. A conjecture of Armstrong, proved recently by Johnson, says that the average size of an $(s,t)$-core, when $s$ and $t$ are coprime positive integers, is $\frac1{24}(s-1)(t-1)(s+t-1)$. Armstrong also conjectured that the same formula gives the average size of a self-conjugate $(s,t)$-core; this was proved by Chen, Huang and Wang.

In the present paper, we develop the ideas from the author's paper \cite{mfcores}, studying actions of affine symmetric groups on the set of $s$-cores in order to give variants of Armstrong's conjectures in which each $(s,t)$-core is weighted by the reciprocal of the order of its stabiliser under a certain group action. Informally, this weighted average gives the expected size of the $t$-core of a random $s$-core.
\end{abstract}

\pdfpagewidth=8.27in
\pdfpageheight=11.69in

\renewcommand\baselinestretch{1.01}

\section{Introduction}

The study of integer partitions is a very active subject, with connections to representation theory, number theory and symmetric function theory. A particularly prominent theme is the study of \emph{$s$-core partitions}, when $s$ is a natural number: we say that a partition $\la$ is an $s$-core if it does not have a rim hook of length $s$; if $\la$ is any partition, then the $s$-core of $\la$ is the partition obtained by repeatedly removing rim $s$-hooks. The set of all $s$-cores displays a geometric structure, with connections to Lie theory. In the case where $s$ is a prime, $s$-cores play an important role in the $s$-modular representation theory of the symmetric group.

In the last few years, there has been considerable interest in the study of $(s,t)$-cores, i.e.\ partitions which are both $s$- and $t$-cores, for given natural numbers $s$ and $t$. When $s$ and $t$ are coprime, there are only finitely many $(s,t)$-cores; the exact number was computed by Anderson \cite{and}, and in the particular case where $t=s+1$ coincides with the $s$th Catalan number. The properties of $(s,t)$-cores have been studied from a variety of aspects: Fishel and Vazirani \cite{fv} explored the connection with alcove geometry and the Shi arrangement, and several authors \cite{kane,os,vand,mfcores} have studied the properties of the unique largest $(s,t)$-core. The present author \cite{mfcores} defined a \emph{level $t$ action} of the affine symmetric group $\weyl s$ on the set of $s$-cores (generalising an action due to Lascoux \cite{l} in the case $t=1$) and showed that two $s$-cores have the same $t$-core \iff they lie in the same orbit for this action.

Recently, Armstrong has examined the sizes of $(s,t)$-cores, conjecturing in \cite{ahj} that the average size of an $(s,t)$-core is given by $\frac1{24}(s-1)(t-1)(s+t-1)$; he made the same conjecture for $(s,t)$-cores which are \emph{self-conjugate}, i.e.\ symmetric down the diagonal. The conjecture for self-conjugate $(s,t)$-cores was proved soon afterwards by Chen, Huang and Wang \cite{chw}, but the original conjecture proved more difficult. The `Catalan case' $t=s+1$ was proved by Stanley and Zanello \cite{sz}, and this was generalised to the case $t\equiv1\ppmod s$ by Aggarwal \cite{agg}. Very recently, the full conjecture was proved by Johnson \cite{j} using Ehrhart theory.

In this paper, we connect Armstrong's conjectures to the level $t$ action of the affine symmetric group on the set of $s$-cores, and present variants of these conjectures, in which the size of an $(s,t)$-core $\la$ is weighted by the reciprocal of the order of the stabiliser of $\la$ under this action. Surprisingly, these weighted averages are (apparently) given by simple formul\ae{} which are very similar to those in Armstrong's conjectures. We motivate our conjectures in terms of choosing an $s$-core at random and asking for the expected size of its $t$-core.

We now indicate the layout of this paper. In the next section we give some basic definitions and recall Armstrong's conjectures. In \cref{nscsec} we consider actions of the affine symmetric group on the set of $s$-cores and give our variant of Armstrong's conjecture. We show how to compute the stabiliser of an $(s,t)$-core, and connect this to Johnson's geometric approach. We then consider actions on certain finite sets of $s$-cores; this allows a rigorous interpretation of our conjecture in terms of the $t$-core of a randomly chosen $s$-core. Finally, we give (with proof) a formula for the denominator in our weighted average, i.e.\ the sum of the reciprocals of the orders of the stabilisers of the $(s,t)$-cores. In \cref{scsec} we consider self-conjugate cores, introducing an action of the affine hyperoctahedral group on the set of self-conjugate $(s,t)$-cores and giving a weighted variant of Armstrong's conjecture in this case. As in the non-self-conjugate case, we show how to compute the stabiliser of a self-conjugate $(s,t)$-core and explore the connections to Johnson's work, before studying actions on finite sets of self-conjugate $s$-cores.

This paper is mostly self-contained, although several results from the author's previous paper \cite{mfcores} are used. We also use some standard results on Coxeter groups without proof.

We remark that since this paper first appeared in preprint form, our main conjectures have been proved by Wang \cite{wang}.

\section{Armstrong's conjectures}\label{armsec}

\bigskip
\begin{mdframed}[innerleftmargin=3pt,innerrightmargin=3pt,innertopmargin=3pt,innerbottommargin=3pt,roundcorner=5pt,innermargin=-3pt,outermargin=-3pt]
\textit{We assume throughout this paper that $s$ and $t$ are coprime natural numbers with $s\gs2$, and we define $\ci st=\frac12(s-1)(t-1)$ and $u=\inp{s/2}$.}
\end{mdframed}
\medskip

In this paper, a \emph{partition} means a weakly decreasing infinite sequence $\la=(\la_1,\la_2,\dots)$ of non-negative integers such that $\la_i=0$ for large $i$. If $\la$ is a partition, we write $\card\la=\la_1+\la_2+\cdots$, and refer to this as the \emph{size} of $\la$. We write $\la'$ for the \emph{conjugate} partition, defined by $\la'_i=\cardx{\lset{j\gs1}{\la_j\gs i}}$, and we say that $\la$ is \emph{self-conjugate} if $\la=\la'$.

The \emph{Young diagram} of a partition $\la$ is the set
\[
[\la]=\lset{(r,c)\in\bbn^2}{j\ls\la_i}.
\]
If $(r,c)\in[\la]$, then the \emph{$(r,c)$-rim hook} of $\la$ is the set of all $(s,d)\in[\la]$ such that $s\gs r$, $d\gs c$ and $(s+1,d+1)\notin[\la]$. The \emph{$(r,c)$-hook length} of $\la$ is the size of this rim hook, which equals $1+(\la_r-r)+(\la_c'-c)$.

We say that $\la$ is an \emph{$s$-core} if none of the hook lengths of $\la$ equals $s$ (or equivalently if none of them is divisible by $s$), and we let $\cores s$ denote the set of all $s$-cores. We say that $\la$ is an \emph{$(s,t)$-core} if it is both an $s$-core and a $t$-core, i.e.\ it lies in $\cores s\cap\cores t$. If $\la$ is any partition, then the \emph{$t$-core of $\la$} is the $t$-core obtained by repeatedly removing rim hooks of length $t$.

It is an easy exercise to show that (given the assumption that $s$ and $t$ are coprime) there are are only finitely many $(s,t)$-cores. More specifically, we have the following enumerative results.

\begin{thm}\indent
\begin{enumerate}
\vspace{-\topsep}
\item
\cite[Theorems 1 \& 3]{and} The number of $(s,t)$-cores is $\mfrac1{s+1}\mbinom{s+t}s$.
\item
\cite[Theorem 1]{fms} The number of self-conjugate $(s,t)$-cores is $\mbinom{\lfloor s/2\rfloor+\lfloor t/2\rfloor}{\lfloor s/2\rfloor}$.
\end{enumerate}
\end{thm}

$(s,t)$-cores have been intensively studied in the last few years. A very recent result is the following, which was conjectured by Armstrong \cite[Conjecture 1.6]{ahj}.

\begin{thmc}{j}{Theorem 1.7}\label{armnsc}
The average size of an $(s,t)$-core is $\mfrac{(s+t+1)(s-1)(t-1)}{24}$.
\end{thmc}

Armstrong also conjectured the same statement for self-conjugate $(s,t)$-cores. This was proved (rather earlier than \cref{armnsc}) by Chen, Huang and Wang.

\begin{thmcx}{chw}\label{armsc}
The average size of a self-conjugate $(s,t)$-core is $\mfrac{(s+t+1)(s-1)(t-1)}{24}$.
\end{thmcx}

The purpose of this paper is to present conjectured variants of these two statements, in which the sizes of the $(s,t)$-cores are weighted in a meaningful way. In \cref{nscsec} we give a weighted version of \cref{armnsc}, and in \cref{scsec} we do the same for \cref{armsc}.

\section{A weighted version of Armstrong's conjecture for $(s,t)$-cores}\label{nscsec}

\subsection{Action of the affine symmetric group on $s$-cores}\label{weylsec}

The weightings in our variant of Armstrong's conjecture are defined using an action of the affine symmetric group which first appeared in \cite{mfcores}. Let $\weyl s$ denote the affine symmetric group of degree $s$; this can be defined as the set of all permutations $w$ of $\bbz$ satisfying the following conditions:
\begin{enumerate}
\item
$w(m+s)=w(m)+s$ for all $m\in\bbz$;
\item
$w(0)+w(1)+\dots+w(s-1)=\mbinom s2$.
\end{enumerate}
We will say that a function $w:\bbz\to\bbz$ is \emph{\spr} if it satisfies condition (1) above. We remark that if $w:\bbz\to\bbz$ is \spr and $X\subset\bbz$ is any transversal of the congruence classes modulo $s$ in $\bbz$, then $w$ satisfies condition (2) (and hence lies in $\weyl s$) if and only if $\sum_{x\in X}w(x)=\sum_{x\in X}x$.

\medskip
$\weyl s$ has a well-known presentation by generators and relations. Before we give this, we establish some conventions of notation: if $a\in\bbz$, then we write $\sz a$ for the set $a+s\bbz=\lset{a+sm}{m\in\bbz}$. Then $\zsz$ is the set $\lset{\sz a}{a\in\bbz}$, which is an abelian group under addition in the usual way. We let $\bbz$ act additively and multiplicatively on $\zsz$ in the natural way, i.e.\ via $\sz a+b=\sz{a+b}$ and $\sz ab=\sz{ab}$, for $a,b\in\bbz$.

Now for each $i\in\zsz$, let $w_i$ be the element of $\weyl s$ defined by
\[
w_i(m)=
\begin{cases}
m+1&\text{if }m\in i-1\\
m-1&\text{if }m\in i\\
m&\text{otherwise}.
\end{cases}
\]
Then $\weyl s$ is generated by $\lset{w_i}{i\in\zsz}$, subject to defining relations
\begin{alignat*}2
w_i^2&=1\qquad&&\text{for each }i,\\
w_i w_j&=w_j w_i&&\text{if }j\neq i\pm1,\\
w_iw_jw_i&=w_jw_iw_j&\qquad&\text{if }j=i+1\text{ and }s>2.
\end{alignat*}
Thus $\weyl s$ is the affine Coxeter group of type $\tilde A_{s-1}$.

We define the \emph{level $t$ action} of $\weyl s$ on $\bbz$, denoted $w\mapsto\act w$, by
\[
\act w_im=\begin{cases}
m+t&\text{if }m\in(i-1)t-\ci st\\
m-t&\text{if }m\in it-\ci st\\
m&\text{otherwise}
\end{cases}
\]
for $i\in\zsz$, $m\in\bbz$. Note that if $t=1$ then $\ci st=0$, so this is just the natural action of $\weyl s$ on $\bbz$. We remark that the term $-\ci st$ is not really necessary in this section (and does not appear in the definition of the level $t$ action given in \cite{mfcores}); it can be removed with an easy modification of the results below. But using the term $-\ci st$ means that the action works well with self-conjugate partitions, which will be useful in \cref{scsec}.

We can use the level $t$ action of $\weyl s$ on $\bbz$ to describe an action on the set of $s$-cores, by using beta-sets. If $\la$ is a partition, then the \emph{beta-set} of $\la$ is the set
\[
\be\la=\lset{\la_i-i}{i\in\bbn}.
\]
It is easy to check that if $\la$ is a partition and $w\in\weyl s$ then $\act w\be\la$ is also the beta-set of a (unique) partition, so we can define a level $t$ action $w\mapsto\check w$ of $\weyl s$ on the set of partitions by
\[
\calb(\check w\la)=\act w\be\la
\]
for every partition $\la$ and every $w\in\weyl s$.

This action was introduced by the author in \cite{mfgencores}, where it was shown that the action preserves the set of $s$-cores. So we can restrict the level $t$ action to give an action (which we also denote $w\mapsto\check w$) on $\cores s$. In the case $t=1$, this action was introduced by Lascoux \cite{l}.

\begin{eg}
Take $s=3$ and $t=2$, and let $\la=(5,3,1^2)$, which is a $3$-core. We have
\begin{align*}
\be\la&=\{4,1,-2,-3,-5,-6,-7,\dots\},\\
\intertext{so that}
\act w_{\sz0}\be\la&=\{4,1,-1,-2,-4,-5,-7,-8,-9,\dots\},\\
\act w_{\sz1}\be\la&=\{2,-1,-3,-4,-5,\dots\},\\
\act w_{\sz2}\be\la&=\{6,3,0,-3,-5,-6,-7,\dots\},
\end{align*}
and hence
\[
\check w_{\sz0}\la=(5,3,2^2,1^2),\qquad
\check w_{\sz1}\la=(3,1),\qquad
\check w_{\sz2}\la=(7,5,3,1).
\]
Note that these partitions are all $3$-cores.
\end{eg}

Now we can state our conjecture. Given an $s$-core $\la$, write $\stab st\la$ for the stabiliser of $\la$ under the level $t$ action of $\weyl s$ on $\cores s$. Then our conjecture gives the average size of an $(s,t)$-core, but with each $(s,t)$-core $\la$ weighted by the reciprocal of $\stb st\la$, as follows.

\begin{conj}\label{main}
\[
\frac{\displaystyle\sum_{\la\in\simcores st}\frac{\card\la}{\stb st\la}}{\displaystyle\sum_{\la\in\simcores st}\frac{1}{\stb st\la}}=\frac{(s-1)(t^2-1)}{24}.
\]
\end{conj}

\subsection{Motivation: the $t$-core of an $s$-core}\label{motsubs}

Here we recall some results which will give some meaning to the weighted average in \cref{main}. We begin with a result of Olsson.

\begin{thmc}{ol}{Theorem 1}\label{olsson}
If $\la$ is an $s$-core, then the $t$-core of $\la$ is also an $s$-core.
\end{thmc}

Taking the $t$-core of an $s$-core therefore gives a map from the set of $s$-cores to the set of $(s,t)$-cores. The next result says that the fibres of this map are determined by the level $t$ action of $\weyl s$.

\begin{propnc}{mfcores}{Proposition 4.2 \& Corollary 4.5}\label{sameorbit}
Suppose $\la$ and $\mu$ are $s$-cores. Then $\la$ and $\mu$ have the same $t$-core if and only if they lie in the same orbit under the level $t$ action of $\weyl s$ on $\cores s$.
\end{propnc}

We can informally interpret the weighted average in \cref{main} as weighting each $(s,t)$-core $\la$ `by the size of the orbit containing $\la$'. In fact, all the orbits are infinite, so this does not strictly make sense, though we will make it rigorous below by working with finite sets of $s$-cores. Thus, where Armstrong's conjecture addresses the question `given a random $(s,t)$-core, what is its expected size?', our weighted version addresses the question `given a random $s$-core, what is the expected size of its $t$-core?' It is surprising that the apparent answer is so simple and so similar to Armstrong's conjecture.

\medskip
Before looking at finite sets of $s$-cores, we define $s$-sets, and examine a connection to Johnson's geometric proof of Armstrong's conjecture.

\subsection{$s$-sets}\label{ssetsec}

It will be useful to encode an $s$-core as a set of $s$ integers. To do this, we use the fact (first observed by Robinson \cite[2.8]{robin4}) that a partition $\la$ is an $s$-core if and only if for every $b\in\be\la$ we have $b-s\in\be\la$. With this in mind, we define the \emph{$s$-set} of an $s$-core $\la$ to be the set $\sset\la=\lset{a_i}{i\in\zsz}$, where $a_i$ is the smallest integer in $i$ but not in $\be\la$, for each $i$. Another way of saying this (as pointed by Wang \cite{wang}) is $\sset\la=(\be\la+s)\setminus\be\la$. $\sset\la$ is then a set of $s$ integers which are pairwise incongruent modulo $s$, and which sum to $\binom s2$. In general, we refer to any set of $s$ integers with these two properties as an $s$-set; it is shown in \cite{mfcores} that any $s$-set is the $s$-set of a unique $s$-core.

\begin{eg}
Suppose $s=3$ and $\la=(6,4,2,1^2)$. Then
\[
\be\la=\{5,2,-1,-3,-4,-6,-7,-8,\dots\},
\]
so that $\sset\la=\{8,0,-5\}$.

Now consider the $3$-set $\{6,1,-4\}$. We construct the corresponding $3$-core by constructing the beta-set
\[
\{3,0,-3,\dots\}\cup\{-2,-5,-8,\dots\}\cup\{-7,-10,-13,\dots\}=\{3,0,-2,-3,-5,-6,-7,\dots\}.
\]
This is the beta-set of the partition $(4,2,1^2)$, which is a $3$-core.
\end{eg}

This bijection between $s$-cores and $s$-sets is used in \cite{mfcores} to describe a geometric structure on the set of $s$-cores. Later we will see a different version of this structure which was used by Johnson in the proof of Armstrong's conjecture.

Note that we can describe the level $t$ action of $\weyl s$ on $\cores s$ using $s$-sets: we have
\[
\sset{\check w\la}=\act w\sset\la
\]
for any $\la\in\cores s$ and $w\in\weyl s$. This will allow us to give a formula for $\stb st\la$ in terms of $\sset\la$ below. First we need to examine the level $t$ action of $\weyl s$ in more detail.

\subsection{Basic results on the level $t$ action}

In this section we make some simple observations about the level $t$ action of $\weyl s$ on $\bbz$. The definition given in \cref{weylsec} specifies $\act w_i$ for each $i\in\zsz$, but it is useful to have an explicit expression for $\act w$ when $w$ is any element of $\weyl s$. This is given by the following easy \lcnamecref{expressbarw}.

\begin{lemma}\label{expressbarw}
Suppose $w\in\weyl s$ and $m\in\bbz$, and choose $i\in\bbz$ such that $it\equiv m+\ci st\ppmod s$. Then $\act w(m)=m+t(w(i)-i)$.
\end{lemma}

Hence we can explicitly determine the image of the level $t$ action.

\begin{propn}\label{imageleveltz}
The level $t$ action of $\weyl s$ on $\bbz$ is faithful, and its image is the set
\[
\lset{x\in\weyl s}{x(m)\equiv m\ppmod t\text{ for all }m\in\bbz}.
\]
\end{propn}

\begin{pf}
To show that the level $t$ action is faithful, observe that for $w\in\weyl s$ we have $\act w(it-\ci st)=tw(i)-\ci st$ for all $i\in\bbz$, by \cref{expressbarw}. Hence if $\act w$ is the identity permutation, then so is $w$.

Now we consider the image of the level $t$ action. Take $w\in\weyl s$; then it is clear from \cref{expressbarw} that $\act w(m)\equiv m\ppmod t$ for all $m$ and that $\act w$ is \spr. For each $m\in\{0,\dots,s-1\}$ let $i_m$ be the element of $\{0,\dots,s-1\}$ such that $i_mt\equiv m+\ci st\ppmod s$. Since $s$ and $t$ are coprime, the map $m\mapsto i_m$ is bijective. So by \cref{expressbarw}
\[
\sum_{m=0}^{s-1}(\act w(m)-m)=t\sum_{i=0}^{s-1}(w(i)-i)=0,
\]
which implies that $\act w\in\weyl s$.

Conversely, suppose $x\in\weyl s$ with $x(m)\equiv m\ppmod t$ for all $m$; then we must show that there is $w\in\weyl s$ with $x=\act w$. For each $i\in\bbz$, let $m=it-\ci st$, and set $w(i)=i+(x(m)-m)/t$. The assumption that $t$ divides $x(m)-m$ for every $m$ means that the values of $w$ lie in $\bbz$. Clearly $w$ is \spr since $x$ is, so to show that $w$ is a permutation of $\bbz$ it suffices to show that if $i,j\in\bbz$ with $i\nequiv j\ppmod s$ then $w(i)\nequiv w(j)\ppmod s$: setting $m=it-\ci st$ and $n=jt-\ci st$, we have $m\nequiv n\ppmod s$ (since $s$ and $t$ are coprime) and hence
\[
t(w(i)-w(j))=(it+x(m)-m)-(jt+x(n)-n)=x(m)-x(n)\nequiv0\pmod s.
\]
Hence $w(i)\nequiv w(j)\ppmod s$.

So $w$ is an \spr permutation of $\bbz$; one can show (by essentially a reverse of the argument in the first part of the proof) that $w(0)+\dots+w(s-1)=\mbinom s2$, so $w\in\weyl s$. By construction $x=\act w$, and we are done.
\end{pf}

\begin{cory}\label{xsetbijec}
Suppose $X$ and $Y$ are $s$-sets and $\phi:X\to Y$ is a bijection such that $\phi(x)\equiv x\ppmod t$ for all $x\in X$. Then there is a unique $w\in\weyl s$ such that $\phi(x)=\act w(x)$ for all $x\in X$.
\end{cory}

\begin{pf}
Since $X$ contains exactly one integer in each congruence class modulo $s$, there is a unique \spr function $v:\bbz\to\bbz$ such that $v|_X=\phi$; this function satisfies $v(m)\equiv m\ppmod t$ for every $m$ since $\phi$ does, so by \cref{imageleveltz} it suffices to show that $v\in\weyl s$. To see that $v$ is a bijection, it suffices to show that $v(m)\nequiv v(n)\ppmod s$ when $m\nequiv n\ppmod s$; since $v$ is \spr we may as well take $m,n\in X$, in which case the result follows because the elements of $Y$ are pairwise incongruent modulo $s$ and $\phi$ is injective.

So $v$ is an $s$-periodic permutation of $\bbz$. Since in addition
\[
\sum_{x\in X}v(x)=\sum_{y\in Y}y=\mbinom s2=\sum_{x\in X}x
\]
with $X$ a transversal of the congruence classes modulo $s$, we have $v\in\weyl s$.
%
\end{pf}

\subsection{$s$-sets and stabilisers}\label{ssetstabsec}

Now recall that $\stab st\la$ denotes the stabiliser of an $s$-core $\la$ under the level $t$ action of $\weyl s$. The next result shows how to compute $\stab st\la$ from $\sset\la$.

\begin{propn}\label{mainstb}
Suppose $\la\in\cores s$. For each $i\in\ztz$, let $c_i=\card{\sset\la\cap i}$. Then
\[
\stb st\la=\prod_{i\in\ztz}c_i!.
\]
\end{propn}

\begin{pf}
The description of the level $t$ action on $\cores s$ in terms of $s$-sets given in \cref{ssetsec} means that $\stb st\la$ equals the number of elements of $\weyl s$ fixing $\sset\la$ setwise under the level $t$ action on $\bbz$. The case $t=1$ of \cref{xsetbijec} implies that for every permutation $v$ of $\sset\la$ there is a unique element of $\weyl s$ extending $v$. By \cref{imageleveltz}, this element lies in the image of the level $t$ action if and only if it fixes every integer modulo $t$, which happens if and only if $v$ fixes every element of $\sset\la$ modulo $t$. Since the level $t$ action of $\weyl s$ is faithful, different permutations of $\sset\la$ correspond to different elements of $\weyl s$, so the size of the stabiliser is just the number of permutations of $\sset\la$ that fix every element modulo $t$.
\end{pf}

Now we show how to interpret this formula geometrically. We work in Euclidean space $\bbr^s$, with coordinates labelled using the set $\zsz$. Following \cite{mfcores} we define the affine subspace
\[
P^s=\rset{x\in\bbr^s}{\textstyle\sum_{i\in\zsz}x_i=\mbinom s2}.
\]
Given an $s$-core $\la$, define a point $x_\la\in P^s$ by defining $(x_\la)_i$ to be the unique element of $\sset\la\cap i$ for each $i\in\zsz$. The one-to-one correspondence between $s$-cores and $s$-sets then gives
\[
\lset{x_\la}{\la\in\cores s}=\lset{x\in P^s}{x_i\in i\text{ for all }i\in\zsz}.
\]
Note that this set is a lattice (or rather, an affine lattice), which we denote $\La_s$. This lattice was introduced (with different conventions) by Johnson \cite{j}, who calls $\La_s$ the \emph{lattice of $s$-cores}. Johnson's construction is central to his proof of Armstrong's conjecture via Ehrhart theory; indeed, Johnson makes the legitimate claim that his paper `establishes lattice point geometry as a foundation for the study of simultaneous core partitions'.

Note that this construction is different from that in \cite{fv,mfcores}, where an $s$-core $\la$ is represented a point $p_\la$ in the dominant region of $P^s$; this construction yields a bijection between $\cores s$ and the set of dominant alcoves in $P^s$, but does not yield a lattice.

The advantage of Johnson's construction is the easy identification of the set of $(s,t)$-cores as the set of points of $\La_s$ lying inside a certain simplex. Define a hyperplane $H_i$ in $P^s$ for each $i\in\zsz$ by
\[
H_i=\lset{x\in P^s}{x_i-x_{i-t}=t}.
\]
Let $\ssc st$ denote the simplex bounded by the hyperplanes $H_i$; that is,
\[
\ssc st=\lset{x\in P^s}{x_i-x_{i-t}\ls t\text{ for all }i\in\zsz}.
\]
Then we have the following.

\begin{lemmac}{j}{Lemma 3.1}\label{simplex}
Suppose $\la$ is an $s$-core. Then $\la$ is also a $t$-core \iff $x_\la\in\ssc st$.
\end{lemmac}

\begin{eg}
Suppose $s=3$ and $t=4$. We illustrate part of the lattice of $3$-cores in \cref{3latt}, where we label each point $x$ of $\La_s$ by its coordinates $x_{\sz0},x_{\sz1},x_{\sz2}$ and also by the corresponding $3$-core. The three lines drawn are the hyperplanes $H_{\sz0},H_{\sz1},H_{\sz2}$, and the triangle bounded by these three lines is $\ssc34$. The $3$-cores corresponding to points of $\La_3$ inside this triangle are precisely the $(3,4)$-cores.
\begin{figure}[ht]\small
\[
\begin{tikzpicture}[scale=.5,rotate=60]
\clip(120:1)++(240:2)circle(12.5);
\draw(3,0)++(120:1)++(240:-1)++(300:15)--++(120:30);
\draw(-3,0)++(120:1)++(240:5)++(0:20)--++(180:30)++(0:10)++(240:10)--++(60:30);
\three012\varnothing
\three31{-1}{(1)}
\three04{-1}{(2)}
\three3{-2}2{(1^2)}
\three0{-2}5{(3,1)}
\three{-3}42{(2,1^2)}
\three6{-2}{-1}{(4,2)}
\three{-3}15{(3,1^2)}
\three34{-4}{(2^2,1^2)}
\three61{-4}{(4,2,1^2)}
\three{-3}7{-1}{(5,3,1)}
\three3{-5}5{(3,2^2,1^2)}
\three07{-4}{(5,3,1^2)}
\three{-3}{-2}8{(6,4,2)}
\three6{-5}2{(4,2^2,1^2)}
\three{-6}45{(3^2,2^2,1^2)}
\three0{-5}8{(6,4,2,1^2)}
\three{-6}72{(5,3,2^2,1^2)}
\three{-6}18{(6,4,2^2,1^2)}
\draw(-1.5,0)++(120:2.5)++(240:2)node[fill=white]{$H_{\sz1}$};
\draw(120:-.5)++(240:3.5)node[fill=white]{$H_{\sz2}$};
\draw(1.5,0)++(120:3)++(240:-2.5)node[fill=white]{$H_{\sz0}$};
\end{tikzpicture}
\]
\begin{caption}
{$4$-cores inside the lattice of $3$-cores}
\label{3latt}
\end{caption}
\end{figure}
\end{eg}

The lattice of $s$-cores is also relevant to our study of the level $t$ action of $\weyl s$ on $s$-cores. For $j\in\zsz$ let $r_j:P^s\to P^s$ denote the reflection (with respect to the usual inner product on $\bbr^s$) in the hyperplane $H_j$. Then, as is well known in the theory of reflection groups, the group $\lyew s:=\lspan{r_j}{j\in\zsz}$ is isomorphic to $\weyl s$, and an isomorphism $\theta:\weyl s\to\lyew s$ may be given by mapping
\[
w_i\longmapsto r_{it-\ci st}
\]
for each $i\in\zsz$. Moreover, this isomorphism connects the level $t$ action of $\weyl s$ on $\cores s$ to the action of $\lyew s$ on the lattice $\La_s$, via the following \lcnamecref{actsimplex}.

\begin{lemma}\label{actsimplex}
If $\la\in\cores s$ and $w\in\weyl s$ then
\[
x_{\check w\la}=\theta(w)x_\la.
\]
\end{lemma}

\begin{pf}
In the case where $w=w_i$ for $i\in\zsz$, this follows directly from the formula for a reflection in $\bbr^s$ and the definition of the level $t$ action on $\cores s$. The case for arbitrary $w$ then follows from the fact that $\theta$ is a homomorphism.
\end{pf}

With this geometric interpretation of the level $t$ action, we can realise the stabiliser $\stab st\la$ geometrically. First we show that $\stab st\la$ is a parabolic subgroup of $\weyl s$.

\begin{lemma}\label{parab}
Suppose $\la\in\cores s\cap\cores t$, and let $I$ be the set of $i\in\zsz$ such that $x_\la\in H_{it-\ci st}$. Then
\[
\stab st\la=\rspan{w_i}{i\in I}.
\]
\end{lemma}

\begin{pf}
Recall from the proof of \cref{mainstb} the correspondence between $\stab st\la$ and the group of permutations of $\sset\la$ that fix every element modulo $t$. It follows from the proof of \cite[Proposition 4.1]{mfcores} that since $\la$ is an $(s,t)$-core the elements of $\sset\la$ lying in a given congruence class modulo $t$ form an arithmetic progression with common difference $t$, say $a,a+t,\dots,a+rt$. The group of permutations of these integers is generated by the transpositions $(a+(k-1)t,a+kt)$ for $1\ls k\ls r$. But the transposition $(a+(k-1)t,a+kt)$ is simply the restriction to $\sset\la$ of $\act w_i$, where $i\in\zsz$ is such that $a+kt\in it-\ci st$.

So $\stab st\la$ is generated by those $w_i$ for which $\sset\la$ contains integers $m,m-t$ with $m\in it-\ci st$. This is exactly the condition that $x_\la\in H_{it-\ci st}$.
\end{pf}

From \cref{parab,actsimplex} we deduce the following, which enables us to calculate $\stb st\la$ purely geometrically.

\begin{cory}\label{stabhyp}
If $\la$ is an $(s,t)$-core, then $\stab st\la$ is isomorphic to the group generated by $\lset{r_j}{x_\la\in H_j}$.
\end{cory}

\begin{eg}
Continuing from the last example, we see that $x_\varnothing$ does not lie on any of the hyperplanes $H_{\sz0},H_{\sz1},H_{\sz2}$, so $\stab34\varnothing$ is trivial. The $3$-cores $(1)$, $(2)$ and $(1^2)$ each lie on only one of the three hyperplanes, so the stabiliser of each of these $3$-cores has order $2$. $(3,1^2)$ lies on $H_{\sz1}$ and $H_{\sz2}$, so its stabiliser is isomorphic to the group generated by $r_{\sz0}$ and $r_{\sz1}$, which has order $6$.

So the weighted average in \cref{main} is
\[
\frac{\frac01+\frac12+\frac22+\frac22+\frac56}
{\frac11+\frac12+\frac12+\frac12+\frac16}=\frac54=\frac{(3-1)(4^2-1)}{24},
\]
and \cref{main} is verified in the case $(s,t)=(3,4)$.
\end{eg}

\subsection{Actions on finite sets of cores}\label{finitesec}

In this section we define a family of finite sets of $s$-cores on which $\weyl s$ acts. This will enable us to give rigorous meaning to our interpretation of \cref{main} in terms of random cores.

Choose $N\in\bbn$, and let $\rcores sN$ denote the set of all $s$-cores $\la$ such that $k-l<Ns$ for all $k,l\in\sset\la$. Equivalently, these are the $s$-cores $\la$ such that $\la_1+\la'_1\ls(N-1)s$.

We begin by enumerating these cores.

\begin{lemma}\label{countrcores}
$\card\rcnn=N^{s-1}$.
\end{lemma}

\begin{pf}
Choosing an element of $\rcnn$ amounts to choosing an $s$-set whose elements differ by less than $Ns$. Define a \emph{shifted $s$-set} to be a set of $s$ integers with exactly one in each congruence class modulo $s$, and with smallest element $0$. Given an $s$-set $X$, there is a unique shifted $s$-set arising as a translation of $X$, and this shifted $s$-set will be contained in the interval $[0,Ns-1]$ if and only if the elements of $X$ differ by less than $Ns$. Conversely, given a shifted $s$-set $Y$, we have $\sum_{x\in Y}x\equiv\binom s2\ppmod s$, so there is a unique $s$-set arising as a translation of $Y$, i.e.\ the translation by $\frac1s\left(\binom s2-\sum_{x\in Y}x\right)$.

So it suffices to count the shifted $s$-sets contained in $[0,Ns-1]$, and clearly there are $N^{s-1}$ of these: for each $1\ls i\ls s-1$, we choose exactly one of the integers $i,i+s,i+2s,\dots,i+(N-1)s$ to be in the set.
\end{pf}

To define an action of $\weyl s$ on $\rcnn$, we will show that $\rcnn$ is a transversal of the equivalence classes for an equivalence relation on $\cores s$ which is fixed by the action of $\weyl s$ on $\cores s$. Given $\la,\mu\in\cores s$, set $\la\eqn\mu$ if there is a bijection $\phi:\sset\la\to\sset\mu$ such that $\phi(k)\equiv k\ppmod{Ns}$ for all $k\in\sset\la$. Then obviously $\eqn$ is an equivalence relation on $\cores s$, and we have the following two results.

\begin{propn}\label{uniquerep}
Each equivalence class in $\cores s$ under the relation $\eqn$ contains a unique element of~$\rcnn$.
\end{propn}

\begin{pf}
First we show that each equivalence class contains at least one element of $\rcnn$, i.e.\ given $\la\in\cores s$, there is $\nu\in\rcnn$ such that $\la\eqn\nu$. We proceed by induction on $\sum_{m\in\sset\la}m^2$. Supposing $\la\notin\rcnn$, there are $k,l\in\sset\la$ such that $k-l>Ns$. The set $\sset\la\cup\{k-Ns,l+Ns\}\setminus\{k,l\}$ is the $s$-set of an $s$-core $\mu$ with $\mu\eqn\la$, satisfying $\sum_{m\in\sset\mu}m^2<\sum_{m\in\sset\la}m^2$. By induction $\mu\eqn\nu$ for some $\nu\in\rcnn$.

For uniqueness, suppose $\la,\mu\in\rcnn$ with $\la\eqn\mu$; then we must show that $\la=\mu$. Let $\phi:\sset\la\to\sset\mu$ be the bijection such that $\phi(k)\equiv k\ppmod{Ns}$ for all $k$. Since $\sset\la$ lies within an interval of length $Ns$ and so does $\sset\mu$, the only possibility is that for every $k,l\in\sset\la$ with $k>l$, either $\phi(k)-k=\phi(l)-l$ or $\phi(k)-k=\phi(l)-l-Ns$. So there are integers $a,b$ with $a\in\{0,\dots,s-1\}$ such that for $k\in\sset\la$
\[
\phi(k)=
\begin{cases}
k+(b-1)Ns&(\text{if $k$ is one of the $a$ largest elements of $\sset\la$})\\
k+bNs&(\text{otherwise}).
\end{cases}
\]
But this gives
\[
\mbinom s2=\sum_{k\in\sset\mu}k=\sum_{k\in\sset\la}k+(bs-a)Ns=\mbinom s2+(bs-a)Ns.
\]
So $a=bs$, and therefore $a=b=0$. So $\sset\la=\sset\mu$, and hence $\la=\mu$.
\end{pf}

\begin{propn}\label{presrel}
Suppose $N\in\bbn$. The equivalence relation $\eqn$ on $\cores s$ is preserved by the level $t$ action of~$\weyl s$.
\end{propn}

\begin{pf}
Suppose $\la\eqn\mu$, and let $\phi:\sset\la\to\sset\mu$ be a bijection such that $\phi(k)\equiv k\ppmod{Ns}$ for all $k\in\sset\la$. Since $\sset\la$ contains exactly one integer in each equivalence class modulo $s$, $\phi$ is in fact the unique bijection such that $\phi(k)\equiv k\ppmod s$ for each $k$.

Now take $i\in\zsz$ and consider applying $\check w_i$ to both $\la$ and $\mu$. Let $k$ be the unique element of $\sset\la\cap(it-\ci st)$, and $l$ the unique element of $\sset\la\cap((i-1)t-\ci st)$. Then $\sset{\check w_i\la}=\sset\la\cup\{k-t,l+t\}\setminus\{k,l\}$. The definition of $\phi$ means that $\phi(k)$ is the unique element of $\sset\mu\cap(it-\ci st)$, and similarly for $\phi(l)$. So $\sset{\check w_i\mu}=\sset\mu\cup\{\phi(k)-t,\phi(l)+t\}\setminus\{\phi(k),\phi(l)\}$.

So we can define a bijection $\psi:\sset{\check w_i\la}\to\sset{\check w_i\mu}$ by
\[
\psi(k-t)=\phi(k)-t,\qquad\psi(l+t)=\phi(l)+t,\qquad\psi(m)=\psi(m)\text{ for all other }m,
\]
and we have $\psi(m)\equiv m\ppmod{Ns}$ for all $m$. So $\check w_i\la\eqn\check w_i\mu$.
\end{pf}

\cref{uniquerep,presrel} enable us to define a level $t$ action of $\weyl s$ on $\rcnn$, which (if $N$ is understood) we denote $w\mapsto\hat w$: given $w\in\weyl s$ and $\la\in\rcnn$, we define $\hat w\la$ to be the unique element of $\rcnn$ in the same $\eqn$-class as $\check w\la$.

\begin{eg}
Suppose $s=3$ and $N=4$. Then there are sixteen cores in $\rcores sN$. We illustrate the level $t$ action for $t=1$ and $2$ in \cref{s3t1,s3t2}. In these diagrams, an edge labelled with $i\in\bbz/3\bbz$ represents the action of $\hat w_i$; if there is no edge labelled $i$ meeting a core $\la$, then $\hat w_i\la=\la$.

\begin{figure}[ht]
{\small
\[
\begin{tikzpicture}[scale=2.2]
\draw(0,0)node[inner sep=1pt](311){$(3,1^2)$};
\draw(30:1)node[inner sep=1pt](211){$(2,1^2)$};
\draw(30:1)++(0,1)node[inner sep=1pt](11){$(1^2)$};
\draw(0,2)node[inner sep=1pt](1){$(1)$};
\draw(0,3)node[inner sep=1pt](emp){$\varnothing$};
\draw(150:1)++(0,1)node[inner sep=1pt](2){$(2)$};
\draw(150:1)node[inner sep=1pt](31){$(3,1)$};
\draw(150:1)++(210:1)node[inner sep=1pt](42){$(4,2)$};
\draw(210:2)node[inner sep=1pt](531){$(5,3,1)$};
\draw(210:3)node[inner sep=1pt](642){$(6,4,2)$};
\draw(210:1)++(0,-1)node[inner sep=1pt](5311){$(5,3,1^2)$};
\draw(0,-1)node[inner sep=1pt](4211){$(4,2,1^2)$};
\draw(0,-1)++(330:1)node[inner sep=1pt](42211){$(4,2^2,1^2)$};
\draw(330:2)node[inner sep=1pt](32211){$(3,2^2,1^2)$};
\draw(330:3)node[inner sep=1pt](332211){$(3^2,2^2,1^2)$};
\draw(30:1)++(330:1)node[inner sep=1pt](2211){$(2^2,1^2)$};
\draw(emp)--(1)\lb0--(2)\lb1--(31)\lb2--(311)\lb1--(211)\lb2--(11)\lb1--(1)\lb2;
\draw(211)--(2211)\lb0--(32211)\lb2--(332211)\lb1;
\draw(32211)--(42211)\lb0--(4211)\lb2--(311)\lb0;
\draw(4211)--(5311)\lb1--(531)\lb0--(642)\lb2;
\draw(531)--(42)\lb1--(31)\lb0;
\end{tikzpicture}
\]}
\caption{The level $1$ action of $\weyl3$ on $\rcores34$}
\label{s3t1}
\end{figure}
\begin{figure}[ht]
{\small\[
\begin{tikzpicture}[scale=2.2]
\draw(0,0)node[inner sep=1pt](311){$(3,1^2)$};
\draw(30:1)node[inner sep=1pt](211){$(2,1^2)$};
\draw(30:1)++(0,1)node[inner sep=1pt](11){$(1^2)$};
\draw(0,2)node[inner sep=1pt](1){$(1)$};
\draw(0,3)node[inner sep=1pt](emp){$\varnothing$};
\draw(150:1)++(0,1)node[inner sep=1pt](2){$(2)$};
\draw(150:1)node[inner sep=1pt](31){$(3,1)$};
\draw(150:1)++(210:1)node[inner sep=1pt](42){$(4,2)$};
\draw(210:2)node[inner sep=1pt](531){$(5,3,1)$};
\draw(210:3)node[inner sep=1pt](642){$(6,4,2)$};
\draw(210:1)++(0,-1)node[inner sep=1pt](5311){$(5,3,1^2)$};
\draw(0,-1)node[inner sep=1pt](4211){$(4,2,1^2)$};
\draw(0,-1)++(330:1)node[inner sep=1pt](42211){$(4,2^2,1^2)$};
\draw(330:2)node[inner sep=1pt](32211){$(3,2^2,1^2)$};
\draw(330:3)node[inner sep=1pt](332211){$(3^2,2^2,1^2)$};
\draw(30:1)++(330:1)node[inner sep=1pt](2211){$(2^2,1^2)$};
\draw(30:1)++(0,2)coordinate(emp11);
\draw(30:2)coordinate(112211);
\draw(330:3)++(0,1)coordinate(2211332211);
\draw(330:2)++(0,-1)coordinate(33221142211);
\draw(0,-2)coordinate(422115311);
\draw(210:2)++(0,-1)coordinate(5311642);
\draw(210:3)++(0,1)coordinate(64242);
\draw(150:2)coordinate(422);
\draw(150:1)++(0,2)coordinate(2emp);
\foreach\x in{5}{
\foreach\y in{2}{
\draw(emp)--(11)\lb2;
\draw(11)--(2211)\lb1;
\draw(2211)--(332211)\lb0;
\draw(332211)--(42211)\lb1;
\draw(42211)--(5311)\lb0;
\draw(5311)--(642)\lb2;
\draw(642)--(42)\lb0;
\draw(42)--(2)\lb2;
\draw(2)--(emp)\lb0;
}
}
\draw(1)--(311)\lc0;
\draw(531)--(311)\lc1;
\draw(32211)--(311)\lc2;
\draw(211)--(2)\lc0;
\draw(31)--(11)\lc0;
\draw(211)--(42211)\lc2;
\draw(31)--(5311)\lc1;
\draw(4211)--(42)\lc1;
\draw(4211)--(2211)\lc2;
\end{tikzpicture}
\]
}
\caption{The level $2$ action of $\weyl3$ on $\rcores34$}
\label{s3t2}
\end{figure}
\end{eg}

Now we consider orbits and stabilisers for the level $t$ action on $\rcnn$. For the rest of this section we specialise to the case where $N$ is divisible by $t$, and we write $nt$ instead of $N$. Our aim is to connect the level $t$ action on $\rcn$ to \cref{main} by showing that each orbit contains a unique $(s,t)$-core, and that the size of the orbit containing an $s$-core $\la$ is inversely proportional to $\stb st\la$.

The first step is to compute the kernel of the action. For the next \lcnamecref{kernelrc} we must exclude some cases where $nt$ is very small.

\begin{propn}\label{kernelrc}
Suppose $n\in\bbn$, and assume $nt>1$ and $nst>4$. Then the kernel of the level $t$ action of $\weyl s$ on $\rcn$ is
\[
\kere n:=\lset{w\in\weyl s}{w(m)\equiv m\ppmod{ns}\text{ for all }m\in\bbz}.
\]
\end{propn}

\begin{pf}
Take $w\in\weyl s$, and suppose first that $w(m)\equiv m\ppmod s$ for all $m$. We claim that for any $\la\in\rcn$ we have $\hat w\la=\la$ if and only if $w\in\kere n$. By \cref{expressbarw} we have $\act w(m)\equiv m\ppmod s$ for all $m$, so the unique bijection $\phi:\sset\la\to\sset{\check w\la}$ satisfying $\phi(x)\equiv x\ppmod s$ for all $x\in X$ is just the restriction of $\act w$ to $\sset\la$. If $w\in\kere n$, then by \cref{expressbarw} $\act w(m)\equiv m\ppmod{nst}$ for all $m$, so $\phi(x)\equiv x\ppmod{nst}$ for all $x\in\sset\la$. So $\la\eqnt\check w\la$, and hence $\hat w\la=\la$. On the other hand, if $w\notin\kere n$, choose an integer $m$ such that $w(m)\nequiv m\ppmod{ns}$, and let $x$ be the element of $\sset\la$ congruent to $mt-\ci st$ modulo $s$; then we have $\act w(x)\nequiv x\ppmod{nst}$, so $\phi(x)\nequiv x\ppmod{nst}$, and hence $\check w\la\neqnt\la$, i.e.\ $\hat w\la\neq\la$. So our claim holds, and in particular $w$ lies in the kernel of the level $t$ action on $\rcores s{nt}$ if and only if $w\in\kere n$.

Now suppose instead that there is an integer $m$ such that $w(m)\nequiv m\ppmod s$; letting $x=mt-\ci st$, we have $\act w(x)\nequiv x\ppmod s$. Let $y=\act w(x)+s$; then obviously we have $y\equiv\act wx\ppmod s$ but (by the assumption that $nt>1$) $y\nequiv\act wx\ppmod{nst}$. If $s\gs3$, then since $y\nequiv x\ppmod s$ there is an $s$-set $X$ containing both $x$ and $y$; the unique bijection $\phi:X\to\act w(X)$ satisfying $\phi(z)\equiv z\ppmod s$ for all $z\in X$ must map $y$ to $\act w(x)$, and in particular $\phi(y)\nequiv y\ppmod{nst}$. So there is no bijection from $X$ to $\act w(X)$ fixing every element modulo $nst$. So if $\mu$ is the $s$-core with $s$-set $X$, then $\mu\neqnt\check w\mu$. Hence if $\la$ is the unique element of $\rcores st$ with $\la\eqnt\mu$, then $\la\neqnt\check w\la$, and hence $\la\neq\hat w\la$ So $w$ is not in the kernel of the level $t$ action of $\weyl s$ on $\rcores s{nt}$.

It remains to consider the case $s=2$. Taking $x$ as above, consider $y=1-x$; then $\{x,y\}$ is a $2$-set, so if $\act w(x)\nequiv1-x\ppmod{2nt}$, then we can repeat the argument from the paragraph above. So suppose $\act w(x)\equiv1-x\ppmod{2nt}$; repeating the argument with $x+2$ in place of $x$, we can also assume that $\act w(x+2)\equiv-1-x\ppmod{2nt}$. But now
\[
-1-x\equiv\act w(x+2)=\act w(x)+2\equiv1-x+2\pmod{2nt},
\]
which gives $2nt\ls 4$, contradicting the assumptions on $n$.
\end{pf}

Now we consider $(s,t)$-cores in $\rcores s{nt}$.

\begin{lemma}\label{rstcores}
Suppose $n\gs1$. Then each $(s,t)$-core lies in $\rcores s{nt}$. If $\la\in\rcores s{nt}$, then the orbit containing $\la$ under the level $t$ action of $\weyl s$ on $\rcores s{nt}$ contains a unique $(s,t)$-core, namely the $t$-core of $\la$.
\end{lemma}

\begin{pf}
Suppose $\xi$ is an $(s,t)$-core, and write the elements of $\sset\xi$ in increasing order as $x_1,\dots,x_s$. Then \cite[Propositions 4.2, 4.3]{mfcores} implies that $x_{i+1}-x_i\ls t$ for each $i$, and this implies that $\xi\in\rcores st$, and hence $\xi\in\rcn$.

For the second part of the \lcnamecref{rstcores}, let $O$ be the orbit containing $\la$. Then by \cite[Propositions 4.2, 4.3]{mfcores} there is $w\in\weyl s$ such that $\mu:=\check w\la$ is the $t$-core of $\la$, and in particular $\mu$ is an $(s,t)$-core. Since by the first part of the \lcnamecref{rstcores} $\mu$ lies in $\rcores s{nt}$, we have $\mu=\hat w\la\in O$, so $O$ contains an $(s,t)$-core. For uniqueness, suppose $O$ contains another $(s,t)$-core $\nu$. Then $\nu=\hat v\mu$ for some $v\in\weyl s$, so $\nu\eqnt\check v\mu$. The definitions of $\check v$ and the relation $\eqnt$ now imply that there is a bijection $\phi:\sset\mu\to\sset\nu$ such that $\phi(x)\equiv x\ppmod t$ for each $x$; so by \cite[Proposition 4.1]{mfcores} $\mu$ and $\nu$ have the same $t$-core. Since $\mu$ and $\nu$ are $t$-cores, this means that $\mu=\nu$.
\end{pf}

Now we look at orbit sizes. Recall that $\stab st\la$ denotes the stabiliser of an $s$-core $\la$ under the level $t$ action of $\weyl s$ on $\cores s$.

\begin{lemma}\label{stabsame}
Suppose $\la\in\cores s$. Then $\stab st\la\cap\kere n=\{1\}$.
\end{lemma}

\begin{pf}
Suppose $w\in\kere n$. Then we have $w(m)\equiv m\ppmod s$ for all $m\in\bbz$, and hence $\act w(m)\equiv m\ppmod s$ for every $m\in\bbz$.

If in addition $w\in\stab st\la$, then we have $\act w(\sset\la)=\sset\la$. But the elements of $\sset\la$ are pairwise incongruent modulo $s$, so in fact we must have $\act w(x)=x$ for every $x\in\sset\la$. Since $\sset\la$ is a transversal of the congruence classes modulo $s$, this gives $w(m)=m$ for every integer $m$, so $w=1$.
\end{pf}

Now for $\la\in\rcores s{nt}$ let $\rstb n\la$ denote the stabiliser of $\la$ under the level $t$ action of $\weyl s$ on~$\rcores s{nt}$.

\begin{lemma}\label{rstab}
Suppose $\la\in\rcores s{nt}$. Then $\rstb n\la=\kere n\stab st\la$.
\end{lemma}

\begin{pf}
Clearly both $\kere n$ and $\stab st\la$ lie inside $\rstb n\la$, so $\kere n\stab st\la$ lies inside $\rstb n\la$ as well. For the opposite containment, suppose $w\in\rstb n\la$. Then by definition $\check w\la\eqnt\la$; let $\phi$ denote the bijection $\sset\la\to\sset{\check w\la}$ such that $\phi(x)\equiv x\ppmod{nst}$ for all $x$. Then there is $y\in\weyl s$ such that $\phi$ is just the restriction to $\sset\la$ of $y$. Moreover, we have $y(m)\equiv m\ppmod t$ for all $m\in\bbz$, so by \cref{imageleveltz} there is $v\in\weyl s$ such that $y=\act v$. Following the construction of $v$ given in the proof of \cref{imageleveltz}, we see that since $y(m)\equiv m\ppmod{nst}$ for every $m$, we have $v(m)\equiv m\ppmod{ns}$ for every $m$; that is, $v\in K^{(n)}$.

Now $\act v(\sset\la)=\phi(\sset\la)=\sset{\check w\la}$, so $\check v\la=\check w\la$. So $v^{-1}w\in\stab st\la$, and so $w\in\kere n\stab st\la$.
\end{pf}

The last two results show that if $\la\in\rcn$ then the size of the orbit containing $\la$ is inversely proportional to $\stb st\la$. In fact, we can be more precise, given the following lemma.

\begin{lemma}\label{indexk}
The index of $\kere n$ in $\weyl s$ is $n^{s-1}s!$.
\end{lemma}

The case $n=1$ of this lemma is very well known in Lie theory; it arises from the fact that the affine symmetric group is the semidirect product of the finite symmetric group with its root lattice.

\begin{pf}
We begin with the case $n=1$. Let $H$ denote the setwise stabiliser of $\{0,\dots,s-1\}$ in $\weyl s$. Then clearly $H\cap\kere1=\{1\}$, while $H\kere1=\weyl s$: given $w\in\weyl s$, there is an element $h\in H$ defined by $h(m)\equiv w(m)\ppmod s$ for all $m\in\{0,\dots,s-1\}$, and we have $h^{-1}w\in\kere1$. So $\idx{\weyl s}{\kere1}=\card H$, which is obviously $s!$.

To go from the case $n=1$ to the general case, we just need to show that $\idx{\kere1}{\kere n}=n^{s-1}$. But $\kere1$ is a free abelian group of rank $s-1$, and $\kere n$ consists of the $n$th powers of the elements in this group, which gives the result.
\end{pf}

This yields the following result giving the sizes of level $t$ orbits in $\rcores s{nt}$.

\begin{cory}\label{orbitsize}
Suppose $n\gs1$ and $\la\in\rcores s{nt}$. Then the size of the orbit containing $\la$ under the level $t$ action of $\weyl s$ on $\rcores s{nt}$ is $\mfrac{n^{s-1}s!}{\stb st\la}$.
\end{cory}

\begin{pf}
The cases where $n=t=1$ or $(s,t,n)=(2,1,2)$ are easy to deal with, so we assume that $nt>1$ and $nst>4$, which enables us to use \cref{kernelrc}. Let $\wsn$ denote the image of the level $t$ action of $\weyl s$ on $\rcores s{nt}$. Then $\wsn=\weyl s/\kere n$ by \cref{kernelrc}, so that $\card\wsn=n^{s-1}s!$ by \cref{indexk}. The stabiliser of $\la$ under the action of $\wsn$ is $\mfrac{\rstb n\la}{\kere n}=\mfrac{\kere n\stab st\la}{\kere n}$, by \cref{rstab}. Hence the order of this stabiliser is
\[
\cardx{\frac{\kere n\stab st\la}{\kere n}}=\cardx{\frac{\stab st\la}{\stab st\la\cap\kere n}}=\stb st\la,
\]
so by the Orbit-Stabiliser Theorem the size of the orbit containing $\la$ is $\mfrac{n^{s-1}s!}{\stb st\la}$.
\end{pf}

This result enables us to make precise our informal motivation from \cref{motsubs} concerning random $s$-cores. We now have $\weyl s$ acting on a finite set $\rcn$, and we can select an $s$-core uniformly randomly from this set. By \cref{rstcores} each orbit contains a unique $(s,t)$-core, which is the common $t$-core of all the partitions in this orbit. Hence if $\la\in\cores s\cap\cores t$ then the probability of choosing an $s$-core whose $t$-core is $\la$ is proportional to the size of the orbit containing $\la$, which in turn is inversely proportional to $\stb st\la$. So the left-hand side of \cref{main} gives the expected size of the $t$-core of $\la$.

\subsection{The denominator}\label{denomsec}

Another consequence of the results in \cref{finitesec} is a formula for the denominator appearing in \cref{main}.

\begin{propn}\label{denom}
\[
\sum_{\la\in\simcores st}\frac{1}{\stb st\la}=\frac{t^{s-1}}{s!}.
\]
\end{propn}

\begin{pf}
We specialise the results of \cref{finitesec} to the case $n=1$. By \cref{countrcores}, $\card{\rcores st}=t^{s-1}$, and this is the sum of the sizes of the orbits of $\weyl s$ on $\rcores st$. Each of these orbits contains a unique $(s,t)$-core, so we just sum the result of \cref{orbitsize} over all $(s,t)$-cores $\la$. We obtain
\[
\sum_{\la\in\simcores st}\frac{s!}{\stb st\la}=t^{s-1},
\]
which gives the result.
\end{pf}

\section{A weighted version of Armstrong's conjecture for self-conjugate $(s,t)$-cores}\label{scsec}

Now we consider analogues of the results and conjectures in the previous section for self-conjugate cores. The structure of this section is largely the same as in \cref{nscsec}, though we will be able to be briefer by using some results from that section.

Throughout this section let $\sccores s$ denote the set of all self-conjugate $s$-cores.  Recall that we define $u=\inp{s/2}$.

\subsection{The affine hyperoctahedral group}\label{affhypsec}

We begin by defining a subgroup of $\weyl s$ that fixes $\sccores s$, and which will take the place of $\weyl s$ in this section. For $i\in\zsz$, define $v_i\in\weyl s$ by
\[
v_i=
\begin{cases}
w_i&(\text{if }i=\sz0)\\
w_iw_{-i}&(\text{if }i=\sz l\text{ or }-\sz l\text{ for }1\ls l<u)\\
w_iw_{-i}w_i&(\text{if }i=\sz u\text{ or }-\sz u).
\end{cases}
\]
Note that we have $v_i=v_{-i}$ for each $i\in\zsz$. Now define $\hyp s=\lspan{v_i}{i\in\zsz}$. Then $\hyp s$ is isomorphic to the affine hyperoctahedral group of degree $u$, i.e.\ the affine Coxeter group of type $\tilde C_u$.

It will be helpful to describe $\hyp s$ explicitly in terms of permutations.

\begin{propn}\label{hyponz}
\[
\hyp s=\lset{w\in\weyl s}{w(-1-m)=-1-w(m)\text{ for all }m\in\bbz}.
\]
\end{propn}

\begin{pf}
Let $H$ denote the given subgroup of $\weyl s$. It is easily checked that each $v_i$ satisfies $v_i(-1-m)=-1-v_i(m)$ for all $m$, so $v_i\in H$, and hence $\hyp s\ls H$.

Conversely, suppose $w\in H$, and define
\[
M(w):=\sum_{i=0}^{u-1}(w(i)-i)^2.
\]
We will prove by induction on $M(w)$ that $w\in\hyp s$. In the case $M(w)=0$, we have $w(i)=i$ for $i=0,\dots,u-1$; the fact that $w$ is $s$-periodic and $w(-1-m)=-1-w(m)$ for all $m\in\bbz$ then means that $w(m)=m$ for all $m\in\bbz$, so $w$ is the identity permutation, which lies in $\hyp s$.

For the inductive step, assume $M(w)>0$ and suppose first that there is $a\in\{1,\dots,u-1\}$ for which $w(a)<w(a-1)$. Let $w'=wv_{\sz a}$; then for $i\in\{0,\dots,u-1\}$ we have
\[
w'(i)=
\begin{cases}
w(i-1)&(i=a)\\
w(i+1)&(i=a-1)\\
w(i)&(\text{otherwise}),
\end{cases}
\]
so that $M(w')=M(w)-2w(a-1)+2w(a)<M(w)$. $w'$ lies in $H$, so by induction $w'$ lies in $\hyp s$, and hence so does $w$.

So we may assume that $w(0)<w(1)<\dots<w(u-1)$. Since $M(w)>0$, this means in particular that either $w(0)<0$ or $w(u-1)>u-1$. In the first case, let $w'=wv_{\sz0}$; then we have $w'(i)=w(i)$ for $i=1,\dots,u-1$, while
\[
w'(0)=w(-1)=-1-w(0),
\]
so that $M(w')=M(w)+2w(0)+1<M(w)$, and again we can use the induction hypothesis.

Finally suppose that we are in the case where $w(u-1)>u-1$. Note that if $s$ is odd then in fact $w(u-1)\gs u+1$; this is because when $s$ is odd the conditions on $w$ give $w(u)=u$, so $w(u-1)$ cannot equal $u$. Whether $s$ is even or odd, we let $w'=wv_{\sz u}$. Now we find that $w'(i)=w(i)$ for $i=0,\dots,u-2$, while $w'(u-1)=s-1-w(u-1)$. We obtain
\[
M(w')-M(w)=
\begin{cases}
2u-1-2w(u-1)&(s\text{ even})\\
4u-4w(u-1)&(s\text{ odd})
\end{cases}
\]
which in either case is negative, so again we can apply the inductive hypothesis.
\end{pf}

By restricting the level $t$ action of $\weyl s$ on $\bbz$, we obtain a level $t$ action of $\hyp s$ on $\bbz$. As with $\weyl s$, we can describe the image of this action explicitly.

\begin{propn}\label{imageleveltzsc}
The image of the level $t$ action of $\hyp s$ on $\bbz$ is
\[
\lset{w\in\hyp s}{w(m)\equiv m\ppmod t\text{ for all }m\in\bbz}.
\]
\end{propn}

\begin{pf}
Let $\hyp{s,t}$ denote the image of $\hyp s$ under the level $t$ action of $\hyp s$ on $\bbz$, and let $H$ denote the given subgroup of $\hyp s$. Recall that $\act v$ denotes the image of $v\in\hyp s$ under the level $t$ action on $\bbz$. From \cref{imageleveltz} we know that $\act v_i\in\weyl s$ and $\act v_i(m)\equiv m\ppmod t$ for all $i\in\zsz$ and $m\in\bbz$. It is easy to check that in addition $\act v_i(-1-m)=-1-\act v_i(m)$ for all $m$, so we have $\act v_i\in H$. Hence  $\hyp{s,t}=\lspan{\act v_i}{i\in\zsz}$ is contained in $H$.

For the converse, we follow the proof of \cref{imageleveltz}. Suppose we are given $x\in H$; for $i\in\bbz$ we let $m=it-\ci st$, and set $w(i)=i+\mfrac{x(m)-m}t$. Then (from the proof of \cref{imageleveltz}) $w\in\weyl s$ and $x=\act w$. Moreover, one can easily check that since $x\in\hyp s$ we have $w\in\hyp s$ too. So $x$ lies in~$\hyp{s,t}$. 
\end{pf}

Now we consider the action of $\hyp s$ on $s$-cores. We begin with the following lemma; note that this is where we really require the term $-\ci st$ in the definition of the level $t$ action of $\weyl s$ on $\bbz$.

\begin{lemma}\label{conjugate}
Suppose $\la$ is an $s$-core and $i\in\zsz$. Under the level $t$ action of $\weyl s$ on $\cores s$ we have $(\check w_i\la)'=\check w_{-i}(\la')$.
\end{lemma}

\begin{pf}
It is well-known and easy to prove that
\[
\be{\la'}=\rset{-1-b}{b\in\bbz\setminus\be\la}
\]
for any partition $\la$. The result now follows from the definition of the level $t$ action of $\weyl s$ via beta-sets.
\end{pf}

\begin{eg}
Suppose $s=5$, $t=2$ and $\la=(4,1^2)$. Then
\[
\be\la=\{3,-1,-2,-4,-5,-6,\dots\}.
\]
To apply $\check w_{\sz1}$ to $\la$, we subtract $2$ from all elements of $\be\la$ congruent to $0$ modulo $5$, and add $2$ to all elements congruent to $3$. We obtain
\[
\{5,0,-1,-4,-5,-6,\dots\},
\]
which is the beta-set of $(6,2^2)$. So $\check w_{\sz1}(4,1^2)=(6,2^2)$.

Now consider $\la'=(3,1^3)$. We have
\[
\be{\la'}=\{2,-1,-2,-3,-5,-6,-7,\dots\}.
\]
To apply $\check w_{\sz4}$, we subtract $2$ from all elements congruent to $1$ modulo $5$, and add $2$ to all elements congruent to $4$. We obtain
\[
\{2,1,-2,-3,-4,-5,-7,-8,-9,\dots\}=\be{(3^2,1^4)},
\]
so $\check w_{\sz4}\la'=(3^2,1^4)=(6,2^2)'$, verifying \cref{conjugate} in this case.
\end{eg}

Using \cref{conjugate} and the relations for the generators $w_i$ given in \cref{weylsec}, we deduce the following.

\begin{cory}\label{hypsc}
If $\la\in\sccores s$ and $w\in\hyp s$, then under the level $t$ action of $\weyl s$ on $\cores s$ we have $\check w\la\in\sccores s$.
\end{cory}

In fact, it is not hard to show that (for any value of $t$) $\hyp s$ is the setwise stabiliser of $\sccores s$ under the level $t$ action of $\weyl s$ on $\cores s$.

So we have a level $t$ action of $\hyp s$ on $\sccores s$, which we may also denote $w\mapsto\check w$. Given $\la\in\sccores s$, let $\stabsc st\la$ denote the stabiliser of $\la$ under this action. Now we can state our main conjecture for self-conjugate cores.

\begin{conj}\label{mainsc}
\[
\frac{\displaystyle\sum_{\la\in\simsccores st}\frac{\card\la}{\stbsc st\la}}{\displaystyle\sum_{\la\in\simsccores st}\frac{1}{\stbsc st\la}}=
\begin{cases}
\dfrac{(s-1)(t^2-1)}{24}
&(\text{if $t$ is odd})\\[9pt]
\dfrac{(s-1)(t^2+2)}{24}
&(\text{if $t$ is even}).
\end{cases}
\]
\end{conj}

The rest of this section follows the structure of \cref{nscsec}:  we begin by giving a formula for $\stbsc st\la$, and examining the connection to Johnson's lattice of $s$-cores. We then consider actions of $\hyp s$ on finite sets of self-conjugate $s$-cores, which will enable us to phrase \cref{mainsc} in terms of random self-conjugate $s$-cores, and to give an explicit formula for the denominator in the weighted average in \cref{mainsc}.

We begin by showing that, as in the non-self-conjugate case, the level $t$ orbit containing a self-conjugate $s$-core $\la$ is determined by the $t$-core of $\la$. First we make a definition: say that an $s$-set $X$ is \emph{symmetric} if $s-1-x\in X$ for every $x\in X$.

\begin{lemma}\label{symsset}
Suppose $\la\in\cores s$. Then $\la=\la'$ if and only if $\sset\la$ is symmetric.
\end{lemma}

\begin{pf}
The relationship between $\be\la$ and $\be{\la'}$ given in the proof of \cref{conjugate} yields $\sset{\la'}=\lset{s-1-x}{x\in\sset\la}$ for any $\la\in\cores s$. The result follows.
\end{pf}

\begin{propn}\label{actsc}
If $\la\in\sccores s$, then the $t$-core of $\la$ lies in the same orbit as $\la$ under the level $t$ action of $\hyp s$ on $\sccores s$.
\end{propn}

\begin{pf}
We follow the last part of the proof of \cite[Proposition 4.3]{mfcores}. Let $O$ be the orbit containing $\la$, and let $\nu$ be a partition in this orbit for which the sum $\sum_{k\in\sset\nu}k^2$ is minimised. If we can show that $\nu$ is a $t$-core, then $\nu$ must be the $t$-core of $\la$ (since the level $t$ action of $\weyl s$ on $\cores s$ preserves the $t$-core of an $s$-core).

Suppose for a contradiction that $\nu$ is not a $t$-core. For each $i\in\zsz$ let $k_i$ be the unique element of $\sset\nu\cap i$. By \cref{simplex}, there must be some $j\in\zsz$ such that $k_j>k_{j-t}+t$. By \cref{symsset} $\sset\nu$ is symmetric, so we also have $k_{t-j-1}>k_{-j-1}+t$. Let $i\in\zsz$ be such that $j=it-\ci st$, and consider several cases.
\begin{itemize}
\item
Suppose $i=\sz0$ or $s$ is even and $i=\sz u$. Then $w_i=v_i\in\hyp s$, so $\check w_i\nu$ lies in $O$. Applying $\act w_i$ to $\sset\nu$ amounts to replacing $k_j$ and $k_{j-t}$ with $k_j-t$ and $k_{j-t}+t$. But then
\[
\sum_{k\in\sset{\check w_i\nu}}k^2-\sum_{k\in\sset\nu}k^2=(k_j-t)^2+(k_{j-t}+t)^2-k_j^2-k_{j-t}^2=-2t(k_j-k_{j-t}-t)<0,
\]
contradicting the choice of $\nu$.
\item
Suppose $i=\sz l$ or $-\sz l$, where $1\ls l\ls u-1$, and consider $\check v_i\nu=\check w_i\check w_{-i}\nu$. The conditions on $l$ mean that $j,j-t,t-j-1,-j-1$ are distinct, so applying $\act v_i$ to $\sset\nu$ amounts to replacing $k_j,k_{j-t},k_{t-j-1},k_{-j-1}$ with $k_j-t,k_{j-t}+t,k_{t-j-1}-t,k_{-j-1}+t$. As in the previous case we get $\sum_{k\in\sset{\check v_i\nu}}k^2<\sum_{k\in\sset\nu}k^2$, a contradiction.
\item
Suppose $s$ is odd and $i=\sz u$. Now consider $\check v_i\nu=\check w_{\sz u}\check w_{-\sz u}\check w_{\sz u}\nu$. We now have $j=\sz u$, so that $j=-j-1$. Hence applying $\act v_i$ to $\sset\nu$ amounts to replacing $k_{j+t}$ and $k_{j-t}$ with $k_{j+t}-2t$ and $k_{j-t}+2t$. As in the previous cases we reach a contradiction.
\item
Finally suppose $s$ is odd and $i=-\sz u$. As in the previous case, we can apply $\check v_i$ and reach a contradiction.\qedhere
\end{itemize}
\end{pf}

Hence we get the following analogue of \cref{sameorbit}.

\begin{cory}\label{sameorbitsc}
Suppose $\la,\mu\in\sccores s$. Then $\la$ and $\mu$ have the same $t$-core if and only if they lie in the same orbit under the level $t$ action of $\hyp s$ on $\sccores s$.
\end{cory}

\begin{pf}
If $\la$ and $\mu$ lie in the same level $t$ orbit of $\hyp s$, then they lie in the same level $t$ orbit of $\weyl s$ and so have the same $t$-core by \cref{sameorbit}. The converse follows from \cref{actsc}.
\end{pf}

\subsection{$s$-sets and stabilisers}\label{ssetstabsecsc}

Next we show how to compute $\stbsc st\la$ from the $s$-set for $\la$, when $\la$ is a self-conjugate $s$-core. The method here is the same as in \cref{mainstb}, but the statement is more complicated.

\begin{propn}\label{mainstbsc}
Suppose $\la\in\sccores s$. For $i\in\bbz$, let $c_i=\card{\sset\la\cap(i+t\bbz)}$.
\begin{enumerate}
\item
If $s$ and $t$ are both odd, let $y=(c_u-1)/2$. Then
\[
\stbsc st\la=2^yy!\prod_{i=u+1}^{u+(t-1)/2}c_i!.
\]
\item
If $t$ is even, let $y=(c_u-1)/2$ and $z=c_{u+t/2}/2$. Then
\[
\stbsc st\la=2^{y+z}y!z!\prod_{i=u+1}^{u+(t-2)/2}c_i!.
\]
\item
If $s$ is even, let $y=c_{u+(t-1)/2}/2$. Then
\[
\stbsc st\la=2^yy!\prod_{i=u}^{u+(t-3)/2}c_i!.
\]
\end{enumerate}
\end{propn}

\begin{pf}
We begin exactly as in the proof of \cref{mainstb}. $\stb st\la$ equals the number of elements of $\hyp s$ fixing $\sset\la$ setwise under the level $t$ action on $\bbz$. The fact that $\sset\la$ is a symmetric $s$-set guarantees that if $v$ is a permutation of $\sset\la$ satisfying $v(s-1-i)=s-1-v(i)$ for all $i\in\sset\la$, then there is a unique element of $\hyp s$ extending $v$. By \cref{hyponz}, this element of $\hyp s$ lies in the image of the level $t$ action if and only if if fixes every integer modulo $t$, which happens if and only if $v$ fixes every element of $\sset\la$ modulo $t$. Since the level $t$ action of $\weyl s$ (and hence the level $t$ action of $\hyp s$) on $\bbz$ is faithful, different permutations of $\sset\la$ correspond to different elements of $\weyl s$, so the size of the stabiliser is just the number of permutations $v$ of $\sset\la$ that satisfy $v(s-1-i)=s-1-v(i)$ and $v(i)\equiv i\ppmod t$ for all $i\in\sset\la$. Call such a permutation \emph{good}.

Given $j\in\ztz$, let $\sset\la_j=\sset\la\cap j$. Then any good permutation must permute $\sset\la_j$. If $j\neq s-1-j$, then any permutation of $\sset\la_j$ can occur as the restriction of a good permutation $v$, and the condition that $v(s-1-i)=s-1-v(i)$ for all $i$ then uniquely determines the restriction of $v$ to $\sset\la_{s-1-j}$. On the other hand, if $j=s-1-j$, then the restriction of a good permutation $v$ to $\sset\la_j$ must itself satisfy $v(s-1-i)=s-1-v(i)$ for all $i$; the number of permutations of $\sset\la_j$ achieving this is $2^yy!$, where $y=\lfloor\frac12|\sset\la_j|\rfloor$.

Combining these observations with an analysis of when $j$ equals $s-1-j$ (which depends on the parities of $s$ and $t$) yields the formul\ae{} in the \lcnamecref{mainstbsc}.
\end{pf}

Now as in \cref{ssetstabsec} we connect this result to the lattice of $s$-cores; in the interests of brevity, we omit some of the details here. Recall the affine space
\[
P^s=\rset{x\in\bbr^s}{\textstyle\sum_{i\in\zsz}x_i=\mbinom s2}
\]
and the lattice of $s$-cores
\[
\La_s=\lset{x_\la}{\la\in\cores s}.
\]
Define
\[
Q^s=\rset{x\in P^s}{x_{-1-i}=s-1-x_i\text{ for all }i\in\zsz}.
\]
By \cref{symsset}, an $s$-core $\la$ is self-conjugate \iff $x_\la\in Q^s$. So the set $\lset{x_\la}{\la\in\sccores s}$ is the lattice $\La_s\cap Q^s$, which we call the \emph{lattice of self-conjugate $s$-cores}.

As with $\La_s$, we can identity the set $\lset{x_\la}{\la\in\sccores s\cap\sccores t}$ geometrically. Recall that $H_j$ denotes the hyperplane in $P^s$ defined by the equation $x_j-x_{j-t}=t$, and that $\ssc st$ is the simplex bounded by these hyperplanes. Define $J_j:=H_j\cap Q^s$ and $\stc st:=\ssc st\cap Q^s$. Then $\stc st$ is bounded by the hyperplanes $J_j$ in $Q^s$, and it is immediate from \cref{simplex} that if $\la\in\sccores s$, then $\la$ is a $t$-core if and only if $x_\la$ lies in $\stc st$.

Note that $J_j=J_{t-1-j}$ for each $j$, so there are only $u+1$ distinct hyperplanes $J_j$. Since $Q^s$ is a $u$-dimensional space, this means that $\stc st$ is a simplex in $Q^s$.

\begin{eg}
Consider the case $(s,t)=(4,5)$. In \cref{45sc} we illustrate part of the lattice of self-conjugate $4$-cores, labelling each point with its coordinates $x_{\sz0},x_{\sz1},x_{\sz2},x_{\sz3}$ and with the corresponding $4$-core. The lines drawn are the hyperplanes $J_{\sz0}$, $J_{\sz1}=J_{\sz3}$ and $J_{\sz2}$. The triangle bounded by these lines is $\stc45$, and the points of $\La_s$ it contains are precisely the points $x_\la$ for $\la$ a self-conjugate $(4,5)$-core.
\begin{figure}[htb]\small
\[
\begin{tikzpicture}[scale=.66]
\clip(-10,-5)rectangle(10,11);
\draw(4,-11)--(4,13);
\draw(-12,-7)--(8,13);
\draw(-12,-1)--(16,-1);
\four{0}{1}{2}{3}{$(\varnothing)$}
\four{4}{1}{2}{-1}{$(1)$}
\four{0}{5}{-2}{3}{$(2,1)$}
\four{4}{5}{-2}{-1}{$(2^2)$}
\four{0}{-3}{6}{3}{$(3,1^2)$}
\four{4}{-3}{6}{-1}{$(3,2,1)$}
\four{-4}{1}{2}{7}{$(4,1^3)$}
\four{8}{1}{2}{-5}{$(5,2,1^3)$}
\four{-4}{5}{-2}{7}{$(4,3,2,1)$}
\four{8}{5}{-2}{-5}{$(5,3^2,1^2)$}
\four{-4}{-3}{6}{7}{$(4^2,2^2)$}
\four{8}{-3}{6}{-5}{$(5,4,3,2,1)$}
\four{0}{9}{-6}{3}{$(6,3,2,1^3)$}
\four{4}{9}{-6}{-1}{$(6,3^2,1^3)$}
\four{-4}{9}{-6}{7}{$(6,5,4,3,2,1)$}
\four{-8}{1}{2}{11}{$(8,5,2^3,1^3)$}
\four{8}{9}{-6}{-5}{$(6^2,4^2,2^2)$}
\four{-8}{5}{-2}{11}{$(8,5,4,3,2,1^3)$}
\four{-8}{-3}{6}{11}{$(8,5^2,3^2,1^3)$}
\four{-8}{9}{-6}{11}{$(8,7,6,5,4,3,2,1)$}
\draw(4,3)node[fill=white]{$J_{\sz0}$};
\draw(-2,3)node[fill=white]{$J_{\sz1}$};
\draw(0,-1)node[fill=white]{$J_{\sz2}$};
\end{tikzpicture}
\]
\begin{caption}
{$5$-cores inside the lattice of self-conjugate $4$-cores}
\label{45sc}
\end{caption}
\end{figure}

We illustrate the case $(s,t)=(5,4)$ similarly in \cref{54sc}.
\begin{figure}[htb]\small
\[
\begin{tikzpicture}[scale=.528]
\clip(-7.5,-11.5)rectangle(12.5,13.5);
\draw(0,-11)--(0,13);
\draw(-7,-11)--++(25,25);
\draw(-10,6)--(15,6);
\five{0}{1}{2}{3}{4}{$(\varnothing)$}
\five{5}{1}{2}{3}{-1}{$(1)$}
\five{0}{6}{2}{-2}{4}{$(2,1)$}
\five{5}{6}{2}{-2}{-1}{$(2^2)$}
\five{0}{-4}{2}{8}{4}{$(4,1^3)$}
\five{5}{-4}{2}{8}{-1}{$(4,2,1^2)$}
\five{-5}{1}{2}{3}{9}{$(5,1^4)$}
\five{10}{1}{2}{3}{-6}{$(6,2,1^4)$}
\five{-5}{6}{2}{-2}{9}{$(5,3,2,1^2)$}
\five{10}{6}{2}{-2}{-6}{$(6,3^2,1^3)$}
\five{0}{11}{2}{-7}{4}{$(7,3,2,1^4)$}
\five{-5}{-4}{2}{8}{9}{$(5^2,2^3)$}
\five{5}{11}{2}{-7}{-1}{$(7,3^2,1^4)$}
\five{10}{-4}{2}{8}{-6}{$(6,5,3,2^2,1)$}
\five{0}{-9}{2}{13}{4}{$(9,5,2^3,1^4)$}
\five{-5}{11}{2}{-7}{9}{$(7,6,4,3,2^2,1)$}
\five{5}{-9}{2}{13}{-1}{$(9,5,3,2^2,1^4)$}
\five{10}{11}{2}{-7}{-6}{$(7^2,4^2,2^3)$}
\five{-5}{-9}{2}{13}{9}{$(9,6^2,3^3,1^3)$}
\five{10}{-9}{2}{13}{-6}{$(9,7,6,4,3^2,2,1^2)$}
\five{0}{16}{2}{-12}{4}{$(12,8,4,3,2^4,1^4)$}
\five{5}{16}{2}{-12}{-1}{$(12,8,4^2,2^4,1^4)$}
\draw(2.5,6)node[fill=white]{$J_{\sz1}$};
\draw(0,3.5)node[fill=white]{$J_{\sz4}$};
\draw(2.5,-1.5)node[fill=white]{$J_{\sz0}$};
\end{tikzpicture}
\]
\begin{caption}
{$4$-cores inside the lattice of self-conjugate $5$-cores}
\label{54sc}
\end{caption}
\end{figure}
\end{eg}

Recall that $\weyl s$ acts on $P^s$ via the map $\theta$, under which $w_i$ maps to the reflection $r_{it-\ci st}$. It can be shown (essentially using \cref{imageleveltzsc}) that the stabiliser of $Q^s$ under this action is precisely $\hyp s$. So we have an action (which we will also denote $\theta$) of $\hyp s$ on $Q^s$. If we let $r_j':Q^s\to Q^s$ denote the reflection in the hyperplane $J_j$, then $\theta$ maps $v_i$ to $r'_{it-\ci st}$ for each $i$. We have an analogue of \cref{parab} for self-conjugate $(s,t)$-cores, from which we may deduce the following analogue of \cref{stabhyp}.

\begin{propn}\label{stabhypsc}
If $\la$ is a self-conjugate $(s,t)$-core, then $\stabsc st\la$ is isomorphic to the group generated by $\lset{r'_j}{x_\la\in J_j}$.
\end{propn}

\begin{eg}
Looking again at the last example with $(s,t)=(4,5)$, we see that for $\la\in\sccores4\cap\sccores5$ we have
\[
\stbsc45\la=
\begin{cases}
1&(\la=\varnothing)\\
2&(\la=(1),(2,1),(2^2),(4,1^3))\\
8&(\la=(6,3^2,1^3)).
\end{cases}
\]
So the weighted average in \cref{mainsc} is
\[
\frac{\frac01+\frac12+\frac32+\frac42+\frac72+\frac{15}8}{\frac11+\frac12+\frac12+\frac12+\frac12+\frac18}=3=\frac{(4-1)(5^2-1)}{24},
\]
and \cref{mainsc} holds in this case.

If we take $(s,t)=(5,4)$ instead, then we obtain
\[
\stbsc54\la=
\begin{cases}
2&(\la=\varnothing,(1),(2^2))\\
4&(\la=(2,1))\\
8&(\la=(4,1^3),(6,3^2,1^3)),
\end{cases}
\]
so the weighted average in \cref{mainsc} is now
\[
\frac{\frac02+\frac12+\frac34+\frac42+\frac78+\frac{15}8}{\frac12+\frac12+\frac12+\frac14+\frac18+\frac18}=3=\frac{(5-1)(4^2+2)}{24}.
\]
\end{eg}

\subsection{Actions on finite sets of self-conjugate cores}\label{finitescsec}

Now we consider actions on finite sets of self-conjugate $s$-cores, following the approach in \cref{finitesec}. Given $N\in\bbn$, let $\rsccores sN=\rcores sN\cap\sccores s$; that is, the set of self-conjugate $s$-cores $\la$ such that $k-l<Ns$ for all $k,l\in\sset\la$.

As with $\rcores sN$, we begin by enumerating the elements of $\rsccores sN$.

\begin{lemma}\label{countrsccores}
$\card\rdnn=N^u$.
\end{lemma}

\begin{pf}
Choosing an element of $\rdnn$ amounts to choosing a symmetric $s$-set whose elements differ by less than $Ns$. We define a \emph{shifted doubled symmetric $s$-set} to be a set of $s$ integers all of the same parity and pairwise incongruent modulo $2s$, which is fixed by the map $x\mapsto-x$.

There is an obvious bijection from symmetric $s$-sets to shifted doubled symmetric $s$-sets, which sends an $s$-set $X$ to $\lset{2x-s+1}{x\in X}$. Moreover, the elements of $X$ differ by less than $Ns$ if and only if the elements of the corresponding shifted doubled symmetric $s$-set differ by less than $2Ns$. So it suffices to count the shifted doubled symmetric $s$-sets contained in $[1-Ns,Ns-1]$. So suppose $Y$ is such a set.

Suppose first that $s$ is odd. Then $|Y|$ is odd, so $Y$ must contain $0$. Furthermore, for each $i\in\{2,4,\dots,s-1\}$ $Y$ must contain exactly one integer in $[1-Ns,Ns-1]$ congruent to $i$ modulo $2s$, and must also contain the negative of this integer. So there are $N^{(s-1)/2}$ possibilities for $Y$.

Now suppose $s$ is even. Then the elements of $Y$ must be odd: if the elements of $Y$ are even, then one of them, say $y$, is divisible by $s$; but then $y$ and $-y$ are congruent modulo $2s$ and both lie in $Y$, which is a contradiction unless $y=0$, but this would imply that $|Y|$ is odd, also a contradiction. Now we can see that there are $N^{s/2}$ possibilities for $Y$: for each $i\in\{1,3,\dots,s-1\}$, $Y$ must contain exactly one of the $N$ integers in $[1-Ns,Ns-1]$ congruent to $i$ modulo $2s$, and must also contain the negative of this integer.
\end{pf}

Next we want to show that the level $t$ action of $\weyl s$ on $\rcnn$ restricts to an action of $\hyp s$ on $\rdnn$. To do this, recall the equivalence relation $\eqn$ from \cref{finitesec}. We have the following analogue of \cref{uniquerep}.

\begin{lemma}\label{uniquerepsc}
Each equivalence class in $\sccores s$ under the relation $\eqn$ contains a unique element of $\rdnn$.
\end{lemma}

\begin{pf}
The uniqueness follows immediately from \cref{uniquerep}. For existence, we follow the proof of \cref{uniquerep}. Suppose $\la\in\sccores s$ but $\la\notin\rdnn$; then there are $k,l\in\sset\la$ such that $k-l>Ns$. We may as well take $k,l$ to be the largest and smallest elements of $\sset\la$, and the fact that $\sset\la$ is symmetric then implies that $k+l=s-1$. As in the proof of \cref{uniquerep} we replace $\sset\la$ with $\sset\nu=\sset\la\cup\{k-Ns,l+Ns\}\setminus\{k,l\}$, and the fact that $k+l=s-1$ guarantees that this $s$-set is symmetric, so $\nu$ is self-conjugate. By induction there is $\mu\in\rdnn$ with $\mu\eqn\nu\eqn\la$.
\end{pf}

Now recall that if $w\in\weyl s$ and $\la\in\rcores sN$, then $\hat w\la$ is defined to be the unique element of $\rcores sN$ for which $\hat w\la\eqn\check w\la$.

\begin{propn}\label{scfiniteaction}
The map $w\mapsto\hat w$ restricts to an action of $\hyp s$ on $\rsccores sN$.
\end{propn}

\begin{pf}
We need to show that if $\la$ is self-conjugate and $w\in\hyp s$, then $\hat w\la$ is self-conjugate. By \cref{hypsc} we have $\check w\la\in\sccores s$, and by definition $\hat w\la$ is the unique element of $\rcores sN$ for which $\hat w\la\eqn\check w\la$. But by \cref{uniquerepsc}, the unique such core $\hat w\la$ lies in $\rsccores sN$.
\end{pf}

We will refer to the action in \cref{scfiniteaction} as the level $t$ action of $\hyp s$ on $\rsccores sN$.

\begin{eg}
Suppose $s=5$ and $N=4$. Then there are sixteen cores in $\rsccores sN$. We illustrate the level $t$ action for $t=1$ and $2$ in \cref{s5t1,s5t2}. In these diagrams, an edge labelled $\sz a$ for $a\in\{0,1,2\}$ represents the action of $\hat v_{\sz a}$; if there is no edge labelled $a$ meeting a core $\la$, then $\hat v_{\sz a}\la=\la$.

\begin{figure}[pht]
{\small
\[
\begin{tikzpicture}[scale=-2.1,rotate=22.5]
\foreach\x in{.707}{
\draw(0,0)node[inner sep=1pt](emp){$\varnothing$}
++(0,1)node[inner sep=1pt](1){$(1)$}
++(\x,\x)node[inner sep=1pt](21){$(2,1)$}
++(1,0)node[inner sep=1pt](22){$(2,2)$}
++(0,1)node[inner sep=1pt](4211){$(4,2,1^2)$}
++(-1,0)node[inner sep=1pt](4111){$(4,1^3)$}
++(-\x,\x)node[inner sep=1pt](51111){$(5,1^4)$}
++(0,1)node[inner sep=1pt](621111){$(6,2,1^4)$}
++(\x,\x)node[inner sep=1pt](7321111){$(7,3,2,1^4)$}
++(1,0)node[inner sep=1pt](7331111){$(7,3^2,1^4)$}
++(\x,-\x)node[inner sep=1pt](633111){$(6,3^2,1^3)$}
++(0,-1)node[inner sep=1pt](53211){$(5,3,2,1^2)$}
++(1,0)node[inner sep=1pt](55222){$(5^2,2^3)$}
++(0,1)node[inner sep=1pt](653221){$(6,5,3,2^2,1)$}
++(\x,\x)node[inner sep=1pt](7643221){$(7,6,4,3,2^2,1)$}
++(1,0)node[inner sep=1pt](7744222){$(7^2,4^2,2^3)$};
}
\draw(emp)--(1)\lb0--(21)\lb1--(22)\lb0--(4211)\lb2--(53211)\lb1--(55222)\lb2--(653221)\lb0--(7643221)\lb1--(7744222)\lb0;
\draw(21)--(4111)\lb2--(4211)\lb0;
\draw(53211)--(633111)\lb0--(653221)\lb2;
\draw(4111)--(51111)\lb1--(621111)\lb0--(7321111)\lb1--(7331111)\lb0--(633111)\lb1;
\end{tikzpicture}
\]
}
\caption{The level $1$ action of $\hyp5$ on $\rsccores54$}
\label{s5t1}
\end{figure}
\begin{figure}[pht]
{\small
\[
\begin{tikzpicture}[scale=-2.1,rotate=22.5]
\foreach\x in{.707}{
\draw(0,0)node[inner sep=1pt](emp){$\varnothing$}
++(0,1)node[inner sep=1pt](1){$(1)$}
++(\x,\x)node[inner sep=1pt](21){$(2,1)$}
++(1,0)node[inner sep=1pt](22){$(2,2)$}
++(0,1)node[inner sep=1pt](4211){$(4,2,1^2)$}
++(-1,0)node[inner sep=1pt](4111){$(4,1^3)$}
++(-\x,\x)node[inner sep=1pt](51111){$(5,1^4)$}
++(0,1)node[inner sep=1pt](621111){$(6,2,1^4)$}
++(\x,\x)node[inner sep=1pt](7321111){$(7,3,2,1^4)$}
++(1,0)node[inner sep=1pt](7331111){$(7,3^2,1^4)$}
++(\x,-\x)node[inner sep=1pt](633111){$(6,3^2,1^3)$}
++(0,-1)node[inner sep=1pt](53211){$(5,3,2,1^2)$}
++(1,0)node[inner sep=1pt](55222){$(5^2,2^3)$}
++(0,1)node[inner sep=1pt](653221){$(6,5,3,2^2,1)$}
++(\x,\x)node[inner sep=1pt](7643221){$(7,6,4,3,2^2,1)$}
++(1,0)node[inner sep=1pt](7744222){$(7^2,4^2,2^3)$};
\draw(22)--(53211)\ld2{.293};
\draw(55222)--(4211)\ld2{.6};
\draw(emp)--(22)\lb1--(4211)\lb0--(621111)\lb1;
\draw(55222)--(7744222)\lb1;
\draw(1)--(51111)\lb2--(7331111)\lb1--(7643221)\lb2;
\draw(21)--(4111)\lb0--(633111)\lb1--(653221)\lb0;
\draw(7744222)--(55222)\lb1--(53211)\lb0--(7321111)\lb1;
}
\end{tikzpicture}
\]
}
\caption{The level $2$ action of $\hyp5$ on $\rsccores54$}
\label{s5t2}
\end{figure}
\end{eg}

Now we consider the kernel of the level $t$ action on $\rsccores st$. As in \cref{finitesec}, we now specialise to the case where $N$ is divisible by $t$. Recall that $\kere n$ is the set of elements of $\weyl s$ that fix every integer modulo $ns$.

\begin{propn}\label{sckernelrc}
Suppose $n\in\bbn$ and that $nt>2$. Then the kernel of the level $t$ action of $\hyp s$ on $\rdn$ is $\kerd n:=\hyp s\cap\kere n$.
\end{propn}

\begin{pf}
Take $w\in\hyp s$, and suppose first that $w(m)\equiv m\ppmod s$ for all $m\in\bbz$. As we saw in the proof of \cref{kernelrc}, for any $\la\in\rcn$ we have $\hat w\la=\la$ if and only if $w\in\kere n$. Restricting attention to $\la\in\rdn$, we find that $w$ lies in the kernel of the level $t$ action of $\hyp s$ on $\rdn$ if and only if $w\in\kere n$.

Now suppose there is $m\in\bbz$ such that $w(m)\nequiv m\ppmod s$; note that if $s$ is odd then by \cref{hyponz} $m\nequiv u\ppmod s$. Setting $x=mt-\ci st$, we have $\act w(x)\nequiv x\ppmod s$, and if $s$ is odd then $x\nequiv u\ppmod s$.

If we can find a symmetric $s$-set containing $x$ and an integer $y$ such that $y\equiv\act w(x)\ppmod s$ but $y\nequiv\act w(x)\ppmod{nst}$, then we can proceed as in the proof of \cref{kernelrc}. Since $x\nequiv u\ppmod s$, we can certainly find symmetric $s$-sets containing $x$. The only situation in which we are not free to take $y=\act w(x)+s$ is if $\act w(x)\equiv-1-x\ppmod s$; in this case, if $X$ is a symmetric $s$-set which contains $x$, then it must contain $s-1-x$, so cannot contain any other integer $y\equiv\act w(x)\ppmod s$. So suppose we are in this situation. If we have $s-1-x\nequiv\act w(x)\ppmod{nst}$ then we can proceed as above, so assume that $s-1-x\equiv\act w(x)\ppmod{nst}$. Repeating the argument with $x+s$ in place of $x$, we can also assume that $-1-x\equiv\act w(x+s)\ppmod{nst}$. But this yields $2s\equiv0\ppmod{nst}$, so that $nt\ls2$, contrary to assumption.
\end{pf}

Next we prove an analogue of \cref{rstcores}.

\begin{lemma}\label{rstcoressc}
Suppose $n\gs1$. Then each self-conjugate $(s,t)$-core lies in $\rsccores s{nt}$. If $\la\in\rsccores s{nt}$, then the orbit containing $\la$ under the level $t$ action of $\hyp s$ on $\rsccores s{nt}$ contains a unique $(s,t)$-core, namely the $t$-core of $\la$.
\end{lemma}

\begin{pf}
The first statement is immediate from the corresponding statement in \cref{rstcores}. For the second part, the uniqueness follows from \cref{rstcores}, since the orbit containing $\la$ under the action of $\hyp s$ on $\rsccores s{nt}$ is contained in a level $t$ orbit of $\weyl s$ on $\rcores s{nt}$. The fact that the $t$-core of $\la$ lies in the same orbit as $\la$ follows from \cref{sameorbitsc}.
\end{pf}

Now, as in \cref{nscsec}, we consider the sizes of orbits. Recall that $\stabsc st\la$ denotes the stabiliser of $\la$ under the level $t$ action of $\hyp s$ on $\sccores s$.

\begin{lemma}\label{stabsamesc}
Suppose $\la\in\sccores s$. Then $\stabsc st\la\cap L^{(n)}=\{1\}$.
\end{lemma}

\begin{pf}
This is immediate from \cref{stabsame}, using the fact that the level $t$ action of $\hyp s$ is just the restriction of the level $t$ action of $\weyl s$.
\end{pf}

Now we consider the index of $L^{(n)}$ in $\hyp s$.

\begin{lemma}\label{indexl}
The index of $L^{(n)}$ in $\hyp s$ is $(2n)^uu!$.
\end{lemma}

\begin{pf}
We begin with the case $n=1$. Let $H$ denote the setwise stabiliser of $\{0,\dots,s-1\}$ in $\hyp s$. Then, exactly as in the proof of \cref{indexk} we have $H\cap L^{(1)}=\{1\}$, while $HL^{(1)}=\hyp s$. So $\idx{\hyp s}{L^{(1)}}=|H|$; this equals the number of permutations $h$ of $\{0,\dots,s-1\}$ such that $h(i)+h(s-1-i)=s-1$ for each $i$, which is $2^uu!$.

For the general case, we again copy the proof of \cref{indexk}: $L^{(1)}$ is a free abelian group of rank $u$, and $L^{(n)}$ consists of the $n$th powers of elements in the group, so $\idx{L^{(1)}}{L^{(n)}}=n^u$.
\end{pf}

Now for $\la\in\rsccores s{nt}$ let $\rstbsc n\la$ denote the stabiliser of $\la$ under the level $t$ action of $\hyp s$ on $\rsccores s{nt}$. Then we have the following analogues of \cref{rstab,orbitsize}, which are proved in exactly the same way.

\begin{lemma}\label{rstabsc}
Suppose $\la\in\rsccores s{nt}$. Then $\rstbsc n\la=L^{(n)}\stabsc st\la$.
\end{lemma}

\begin{cory}\label{orbitsizesc}
Suppose $n\gs1$ and $\la\in\rsccores s{nt}$. Then the size of the orbit containing $\la$ under the level $t$ action of $\hyp s$ on $\rsccores s{nt}$ is $(2n)^uu!/\stbsc st\la$.
\end{cory}

As with \cref{rstab}, this result enables us to give a rigorous interpretation of \cref{mainsc} in terms of random self-conjugate $s$-cores. Given $n\gs1$, consider the level $t$ action of $\hyp s$ on $\rdn$. By \cref{rstcoressc} each orbit contains a unique $(s,t)$-core, which is the common $t$-core of all the partitions in this orbit. Hence if we select $\la\in\rsccores s{nt}$ uniformly randomly, then the left-hand side of \cref{mainsc} gives the expected size of the $t$-core of $\la$.

\subsection{The denominator}

As in \cref{denomsec}, we derive a formula for the denominator appearing in \cref{mainsc}. This is proved in exactly the same way as \cref{denom}.

\begin{propn}\label{denomsc}
\[
\sum_{\la\in\simsccores st}\frac{1}{\stbsc st\la}=\frac{t^u}{2^uu!}.
\]
\end{propn}

\end{document}